\edef\savecatcodeat{\the\catcode`@}
\catcode`\@=11

\def\tb@ifSpecChars#1#2{#1}
\def\tb@ifNoSpecChars#1#2{#2}

\def\tableau{%
  \bgroup% matched in \tb@tableauD
  \@ifstar{\let\Tif\tb@ifNoSpecChars\tb@tableauB}% *, don't use special chars
          {\let\Tif\tb@ifSpecChars\tb@tableauB}}% no *, use special chars

\def\tb@tableauB{% add [] if no [options]
  \@ifnextchar[{\tb@tableauC}{\tb@tableauC[]}}

\def\tb@tableauC[#1]{\hbox\bgroup%
    \let\\=\cr% end line
    \def\bl{\global\let\tbcellF\tb@cellNF}%
    \def\tf{\global\let\tbcellF\tb@cellH}% highlighted cell
    \dimen2=\ht\strutbox \advance\dimen2 by\dp\strutbox%
    \ifx\baselinestretch\undefined\relax%
    \else%
       \dimen0=100sp \dimen0=\baselinestretch\dimen0%
       \dimen2=100\dimen2 \divide\dimen2 by\dimen0%
    \fi%
    \let\tpos\tb@vcenter% default position
    \tb@initYoung% default tableau type
    \tb@options#1\eoo% parse options
    \let\arrow\tb@arrow%
    \dimen0=\Tscale\dimen2%
    \dimen1=\dimen0 \advance\dimen1 by \tb@fframe%
    \lineskip=0pt\baselineskip=0pt% line spacing will be from \vbox to \dimen0
    \def\tb@nothing{}%
    \def\endcellno{$\rss\egroup\bss\egroup}% end cell w/o overlap
    \def\endcell{\endcellno\kern-\dimen0}% end cell & prepare to overlap it
    \def\begincell{\vbox to\dimen0\bgroup\vss\hbox to\dimen0\bgroup\hss$}%
    \let\overlay\tb@overlay%
    \let\fl\tb@fl%
    \let\lss\hss\let\rss\hss\let\tss\vss\let\bss\vss% cell alignment
    \def\mkcell##1{% format individual cell
        \let\tbcellF\tb@cellD% default cell frame
        \def\tb@cellarg{##1}% store cell contents
        \ifx\tb@cellarg\tb@nothing\let\tb@cellarg\tb@cellE\fi%
        \begincell\tb@cellarg\endcellno% the actual cell content
        \tbcellF}% draw cell frame
    \let\savecellF\tbcellF% save global value of cellF in case of nested tableau
     \Tif{\catcode`,=4\catcode`|=\active}{}\tb@tableauD}%

\let\tb@savetableauD\tableauD% save any current definition
{% set up characters which will be interpreted as command characters
    \catcode`|=\active \catcode`*=\active \catcode`~=\active%
    \catcode`@=\active% command characters
\gdef\tableauD#1{%
  \Tif{% make all the command characters active in math mode when #1 parsed
    \mathcode`|="8000 \mathcode`*="8000%
    \mathcode`~="8000 \mathcode`@="8000%
    \def@{\bullet}%
    \let|\cr% end line
    \let*\tf% highlighted cell
    \let~\sk% skew cell
  }{}%
  \tpos{\tabskip=0pt\halign{&\mkcell{##}\cr#1\crcr}}%
  \global\let\tbcellF\savecellF% restore global value
  \egroup% match \hbox\bgroup at start of \tableauC
  \egroup}% match \bgroup at start of \tableau
}
\let\tb@tableauD\tableauD% rename the command
\let\tableauD\tb@savetableauD% restore old command with this name
\let\tb@savetableauD\undefined

\def\tb@options#1{\ifx#1\eoo\relax\else\tb@option#1\expandafter\tb@options\fi}

\def\tb@option#1{%
  \if#1t\let\tpos\tb@vtop\fi%        t = align at top
  \if#1c\let\tpos\tb@vcenter\fi%     c = align at center
  \if#1b\let\tpos\vbox\fi%           b = align at bottom
  \if#1F\tb@initFerrers\fi%          F = Ferrers diagram
  \if#1Y\tb@initYoung\fi%            Y = Young diagram
  \if#1s\tb@initSmall\fi%            s = small boxes
  \if#1m\tb@initMedium\fi%           m = medium boxes
  \if#1l\tb@initLarge\fi%            l = large boxes
  \if#1p\tb@initPartition\fi%            p = small partition sized boxes
  \if#1a\tb@initArrow\fi%            a = use arrow font as base dimension
}

\def\tb@vcenter#1{\ifmmode\vcenter{#1}\else$\vcenter{#1}$\fi}

\def\tb@vtop#1{\hbox{\raise\ht\strutbox\hbox{\lower\dimen0\vtop{#1}}}}

\def\tb@initPartition{\def\Tscale{.3}}
\def\tb@initSmall{\def\Tscale{1}}
\def\tb@initMedium{\def\Tscale{2}}
\def\tb@initLarge{\def\Tscale{3}}

\def\tb@initArrow{\dimen2=1.25em}

\def\tb@initYoung{%
  \def\tb@cellE{}% empty cell stays empty
  \let\tb@cellD\tb@cellN% default frame is normal frame
  \def\sk{\global\let\tbcellF\tb@cellNF}}% skew cells are empty
\def\tb@initFerrers{%
  \def\tb@cellE{\bullet}% empty cell gets bullet
  \let\tb@cellD\tb@cellNF% default frame is no frame
  \def\sk{\bullet}}% skew cell gets bullet

\tb@initMedium% default scale

\def\tb@sframe#1{%
  \vbox to0pt{%            Embed frame in a box of no vert or hor extent
    \vss%                            pull box above reference point
    \hbox to0pt{%
      \hss%                          pull box left of reference point
      \vbox to\dimen1{%              Actual width of frame
        \hrule depth #1 height0pt% draw top edge of frame
        \vss%                     begin vcenter sides
        \hbox to\dimen1{%           horiz box with side edges just inside
          \vrule width #1 height\dimen1% left edge
          \hss%                     stretch center
          \vrule width #1%         right edge
          }%
        \vss%                     end vcenter sides
        \hrule height #1 depth 0in% bottom edge
        }%
      \kern-\tb@hframe%           horiz alignment off by half line width
      }%
    \kern-\tb@hframe}}%           vert alignment off by half line width
\def\tb@hframe{.2pt}\def\tb@fframe{.4pt}\def\tb@bframe{2pt}
\def\tb@cellH{\tb@sframe{\tb@bframe}}       % bold frame
\def\tb@cellNF{}                            % no frame
\def\tb@cellN{\tb@sframe{\tb@fframe}}       % normal frame
\let\tbcellF\tb@cellN                       % default is normal

\def\tb@rpad{1pt}
\def\tb@lpad{1pt}
\def\tb@tpad{1.8pt}
\def\tb@bpad{1.8pt}

\def\tb@overlay{\endcell\@ifnextchar[{\tb@overlaya}{\begincell}}
\def\tb@overlaya[#1]{\vbox to\dimen0\bgroup%
  \tb@overlayoptions#1\eoo%
  \tss\hbox to\dimen0\bgroup\lss$}
\def\tb@overlayoptions#1{\ifx#1\eoo\relax\else\tb@overlayoption#1\expandafter\tb@overlayoptions\fi}

\def\tb@overlayoption#1{
  \if#1t\def\tss{\vskip\tb@tpad}\let\bss\vss\fi% t = align at top
  \if#1c\let\tss\vss\let\bss\vss\fi%             c = align at center
  \if#1b\def\bss{\vskip\tb@bpad}\let\tss\vss\fi% b = align at bottom
  \if#1l\def\lss{\hskip\tb@lpad}\let\rss\hss\fi% l = align at left
  \if#1m\let\lss\hss\let\rss\hss\fi%             m = align at middle
  \if#1r\def\rss{\hskip\tb@rpad}\let\lss\hss\fi% r = align at right
}

\def\tb@fl{\endcell\begincell\vrule depth 0pt width \dimen0 height \dimen0 \endcell\begincell}

\@ifundefined{diagram}{}{
\def\tb@arrowpad{.5}

\newoptcommand{\tb@arrow}{\@ne}[2]{%
  \endcell% end previous cell contents
   \begingroup%
   \let\dg@getnodesize\tb@getnodesize% substitute routine to get nodesize
   \dg@USERSIZE=#1\relax%
   \ifnum\dg@USERSIZE<\@ne \dg@USERSIZE=\@ne \fi%
   \dg@parse{#2}%
   \dg@label{\tb@draw{#1}{#2}}}% draw arrow

\def\tb@getnodesize#1#2#3#4#5{\dimen3=\tb@arrowpad\dimen2 #4=\dimen3 #5=\dimen3\relax}
\def\tb@getnodesize#1#2#3#4#5{\ifnum#2=0\ifnum#3=0\tb@getnodesizetail{#4}{#5}\else\tb@getnodesizehead{#4}{#5}\fi\else\tb@getnodesizehead{#4}{#5}\fi}
\def\tb@getnodesizetail#1#2{\dimen3=.5\dimen2 #1=\dimen3 #2=\dimen3}
\def\tb@getnodesizehead#1#2{\dimen3=.5\dimen2 #1=\dimen3 #2=\dimen3}

\def\tb@draw#1#2#3#4{%
        \dg@X=0\dg@Y=0\dg@XGRID=1\dg@YGRID=1\unitlength=.001\dimen0%
        \dg@LBLOFF=\dgLABELOFFSET \divide\dg@LBLOFF\unitlength%
        \dg@drawcalc% compute arrow geometry
        \begincell% start tableau cell
        \let\lams@arrow\tb@lams@arrow% substitute routine
        \begin{picture}(0,0)\begingroup\dg@draw{#1}{#2}{#3}{#4}\end{picture}%
        \endcell% end tableau cell
        \endgroup% match \begingroup in \tb@arrow
        \begincell}% start new entry in this cell
}

\def\tb@lams@arrow#1#2{%
 \lams@firstx\z@\lams@firsty\z@
 \lams@lastx#1\relax\lams@lasty#2\relax
 \lams@center\z@
 \N@false\E@false\H@false\V@false
 \ifdim\lams@lastx>\z@\E@true\fi
 \ifdim\lams@lastx=\z@\V@true\fi
 \ifdim\lams@lasty>\z@\N@true\fi
 \ifdim\lams@lasty=\z@\H@true\fi
 \NESW@false
 \ifN@\ifE@\NESW@true\fi\else\ifE@\else\NESW@true\fi\fi
 \ifH@\else\ifV@\else
  \lams@slope
  \ifnum\lams@tani>\lams@tanii
   \lams@ht\ten@\p@\lams@wd\ten@\p@
   \multiply\lams@wd\lams@tanii\divide\lams@wd\lams@tani
  \else
   \lams@wd\ten@\p@\lams@ht\ten@\p@
   \divide\lams@ht\lams@tanii\multiply\lams@ht\lams@tani
  \fi
 \fi\fi
 \ifH@  \lams@harrow
 \else\ifV@ \lams@varrow
 \else \lams@darrow
 \fi\fi
}

\catcode`\@=\savecatcodeat
\let\savecatcodeat\undefined

\documentclass[tbtags,reqno]{amsart}        

\usepackage[tbtags]{amsmath}

\usepackage{amsthm,amsopn,amsfonts,amssymb,epsfig,pictex}            

\numberwithin{equation}{section}

\newtheorem{theorem}{Theorem}
\newtheorem{lemma}[theorem]{Lemma}
\newtheorem{proposition}[theorem]{Proposition}

\newtheorem{corollary}[theorem]{Corollary}
 
\newtheorem{definition}[theorem]{Definition}

\newtheorem{identity}[theorem]{Identity}
 
\newtheorem{property}[theorem]{Property}

\newtheorem{formula}[theorem]{Formula}
 
\theoremstyle{remark}

\newtheorem*{acknow}{\bf Acknowledgments}

\def\P{{\mathcal P}}

\def\endprf {\square}

\begin{document}

\title[Schur function identities, their $t$-analogs, and 
$k$-Schur irreducibility]
{Schur function identities, their $t$-analogs, and $k$-Schur 
irreducibility }

\author{L. Lapointe}
\address
{McGill University, 
Department of Mathematics and Statistics,  Montr\'eal, Qu\'ebec,
Canada, H3A 2K6}
\email{lapointe@math.mcgill.ca}

\author{J. Morse}
\thanks{Research supported in part by NSF grant \#0100179}
\address
{University of Pennsylvania,
Department of Mathematics,
Philadelphia, PA 19104}
\email{morsej@math.upenn.edu}
 
\begin{abstract}  
We obtain general identities for the product of two Schur
functions in the case where one of the functions is indexed by a 
rectangular partition, and give their $t$-analogs 
using vertex operators.  We study subspaces forming a filtration 
for the symmetric function space that lends itself to generalizing 
the theory of Schur functions and also provides a convenient environment 
for studying the Macdonald polynomials.  We use our identities
to prove that the vertex operators leave such subspaces invariant.
We finish by showing that these operators act simply on the $k$-Schur 
functions, thus leading to a concept of irreducibility for these functions.
\end{abstract}

\maketitle

\section{Introduction }

Let $\Lambda$ be the ring of symmetric functions
in the variables $x_1,x_2,\ldots$, with coefficients
in $\mathbb Q(q,t)$.
The Schur functions, $s_\lambda[X]$,  form a fundamental basis of
$\Lambda$, with central roles in fields such as representation
theory and algebraic geometry.  For example, the Schur functions 
can be identified with the characters of irreducible representations
of the symmetric group, and their products are equivalent to the
Pieri formulas for multiplying Schubert varieties in the intersection ring
of a Grassmannian.  Furthermore, the connection coefficients of
the Schur function basis with various bases
such as the homogeneous symmetric functions,
the Hall-Littlewood polynomials, and the Macdonald polynomials,
are positive and have representation
theoretic interpretations.
In the case of the Macdonald polynomials, $H_\lambda[X;q,t]$,
this expansion takes the form
\begin{equation} \label{MacKos}
H_{\lambda}[X;q,t] = \sum_{\mu} K_{\mu \lambda}(q,t) \, s_{\mu}[X] \, , \qquad
K_{\mu \lambda}(q,t) \in \mathbb N[q,t] \, ,
\end{equation}
where $K_{\mu \lambda}(q,t)$ are known as the $q,t$-Kostka polynomials. 
The representation theoretic interpretation for these polynomials
is given in \cite{[Ga],[Ha]}. 

Recent developments \cite{[LLM],[LM2]} have suggested
that a certain filtration of $\Lambda$ provides
a convenient environment for the generalization 
of natural properties held by the Schur functions, and for 
the study of the Macdonald polynomials.
This filtration, $\Lambda^{(1)}_t \subseteq \Lambda^{(2)}_t \subseteq \cdots 
\subseteq \Lambda^{(\infty)}_t=\Lambda$, is given by 
\begin{equation} 
\Lambda^{(k)}_t = 
{\mathcal L} \left\{ H_\lambda[X;t] \right\}_{\lambda;\lambda_1\leq k} 
=
{\mathcal L} \left\{ H_\lambda[X;q,t] \right\}_{\lambda;\lambda_1\leq k} \, ,
\end{equation}
where $H_\lambda[X;t]$ denote the Hall-Littlewood polynomials.
Two bases for these spaces are introduced in \cite{[LM2]};
the $k$-split polynomial basis,
which is related to a $t$-generalization of Schur function
products \cite{[Shi],[S2],[S1]}, and the $k$-Schur function basis, 
$s_{\lambda}^{(k)}[X;t]$ \cite{[LLM],[LM2]}.  The latter basis plays a 
role in $\Lambda_t^{(k)}$ analogous to the one played in 
$\Lambda$ by the Schur functions.  For example,
work related to the $k$-Schur functions 
prompted a $k$-analog of partition conjugation, a refinement of 
formula (\ref{MacKos}), 
\begin{equation}
H_{\lambda}[X;q,t] = \sum_{\mu} K^{(k)}_{\mu \lambda}(q,t) \, 
s^{(k)}_{\mu}[X;t] \, ,
\end{equation}
where $K^{(k)}_{\mu \lambda}(q,t)$ are conjectured to be in
$\mathbb N[q,t]$,
and a generalization of the Pieri and Littlewood-Richardson rules.

Here, we study classical Schur function properties (and their $t$-analogs)
in the context of the spaces $\Lambda_t^{(k)}$.  In particular, 
the expansion of a product of two Schur functions in terms of Schur 
functions is explicitly known.  In our case, we examine the 
expansion of such a product, where one Schur 
function is indexed by a rectangular partition, in terms of 
certain products of Schur functions.
We find, for nonnegative integers $a,r$ and $m$
and partition $\nu$ with $\ell(\nu)\leq r$, 
\begin{equation}
s_{a^{r+m}} \,  s_{\nu}
 = \sum_{\mu: \, \ell(\mu) \leq r , \, \mu_1 \leq m} 
(-1)^{|\mu|} s_{a^r+\mu} \, 
s_{a-\mu'_m,\ldots,a-\mu'_1,\nu}
\, ,
\end{equation}
where the summand vanishes if $(a-\mu'_m,\ldots,a-\mu'_1)$ is 
not a partition.
  We also generalize this to
 \begin{equation}
B_{a^{r+m}} \,  B_{\nu}
 = \sum_{\mu: \, \ell(\mu) \leq r , \, \mu_1 \leq m} 
(-t)^{|\mu|} B_{a^r+\mu} \, B_{a-\mu'_m,\ldots,a-\mu'_1,\nu}
\, ,
\end{equation}
where $B_{\lambda}$ is a vertex operator that reduces to $s_{\lambda}$
when $t=1$ \cite{[SZ]}.

This result enables us to provide a thorough analysis of 
the operators $B_\lambda$ in the context of our filtration.
We derive commutation relations on $B_\lambda$ as well as a 
number of other identities for these operators and the
Schur functions.  We find that operators (and
Schur functions) indexed by rectangular partitions,
i.e.\! partitions of the form $(\ell^{k+1-\ell})$,
play a particularly important role in our study.
For example, these operators leave fundamental 
subspaces of $\Lambda_t^{(k)}$ invariant.

Results concerning the $k$-Schur function basis
arise as a consequence of our work with operators 
indexed by rectangular partitions.
In the last section, we prove that the action of such
an operator on a $k$-Schur function
produces only one $k$-Schur function.  Namely,
for $\ell=1,2,\ldots,k$,
\begin{equation}
B_{\ell^{k+1-\ell}} \, s_{\lambda}^{(k)}[X;t] = t^d s_\mu^{(k)}[X;t] \, ,
\end{equation}
where $\mu$ is the partition rearrangement of
the entries in $(\ell^{k+1-\ell})$ and $\lambda$, and $t^d$ 
is a positive power of $t$ given explicitly in Theorem~\ref{theorec}.
This result has the important consequence of simplifying the construction of the 
$k$-Schur functions.  In effect, for each $k$, there is a subset 
of $k!$ $k$-Schur functions called the irreducible $k$-Schur functions, 
from which all other $s_\lambda^{(k)}[X;t]$ may be constructed by 
successive application of operators indexed by rectangular partitions.
That is,
\begin{equation}
\label{atirr}
s_{\lambda}^{(k)}[X;t] = t^c B_{R_1} \cdots B_{R_j}\, s_{\mu}^{(k)}[X;t] \, ,
\qquad c \in \mathbb N \, , 
\end{equation}
where $s_{\mu}^{(k)}[X;t]$ is an irreducible $k$-Schur function
and $R_1,\dots,R_{j}$ are rectangular partitions.

Since the Hall-Littlewood polynomials at $t=1$ are the
homogeneous symmetric functions, $h_{\lambda}[X]$,
$\Lambda^{(k)}_t$ reduces to the polynomial
ring $\Lambda^{(k)}=\mathbb Q[h_1,\ldots,h_k]$.
Since $B_R$ is simply multiplication by the Schur function $s_R$
when $t=1$, relation \eqref{atirr} reduces to
\begin{equation}
s^{(k)}_\lambda[X]
=
s_{R_1}[X]
s_{R_2}[X]
\ldots
s_{R_\ell}[X]
s^{(k)}_\mu[X]
\,.
\label{surmul}
\end{equation}
It follows that the irreducible $k$-Schur functions
thus constitute a natural basis for the quotient ring  
$\Lambda^{(k)}_t/{\mathcal I}_k$, where
${\mathcal I}_k$ is the ideal generated by
Schur functions indexed by rectangular shapes of the type $(\ell^{k+1-\ell})$.

\begin{acknow}
{\it 
The enthusiasm from A. Garsia and A. Lascoux greatly contributed 
to this work and we are thankful to M. Zabrocki
for helping us with \cite{[S1]}.  L. Lapointe thanks 
L. Vinet for his support.  J. Morse held an NSF grant
for part of the period devoted to this research.
ACE \cite{[V]} was instrumental towards this work.
}
\end{acknow}

\section{Definitions} \label{sec1}

\subsection{Partitions}
Symmetric polynomials are indexed by partitions, sequences of
non-negative integers $\lambda =(\lambda_1,\lambda_2,\ldots)$ with
$\lambda_1 \ge \lambda_2 \ge \dots$.  
The number of non-zero parts in $\lambda$ is denoted $\ell(\lambda)$
and the degree of $\lambda$ is $|\lambda| = \lambda_1 + \dots + 
\lambda_{\ell(\lambda)}$.  We use $\lambda_L$ to denote 
$\lambda_{\ell(\lambda)}$.  ${\P}^{r}_{\leq m}$ denotes the set of 
all partitions of length at most $r$, and whose first part is not 
larger than $m$.  In this fashion, ${\P}^{r}_{m}$ is the set of partitions 
of length at most $r$, and whose first part is equal to $m$.  The case 
${\P}^{r}_{\leq \infty}$ will be denoted ${\P}^{r}$.
Finally, for a partition $\mu=(\mu_1,\dots,\mu_m)$ of $m$ entries 
(possibly zero), the reverse
reading is denoted $\mu^R=(\mu_m,\dots,\mu_1)$.

We use the dominance order on partitions with $|\lambda|=|\mu|$, where
$\lambda\leq\mu$ when $\lambda_1+\cdots+\lambda_i\leq
\mu_1+\cdots+\mu_i$ for all $i$.  
Given two partitions $\lambda$ and $\mu$,
$\lambda \cup \mu$ stands for the partition 
rearrangement of the parts of $\lambda$ and $\mu$, 
$\lambda\pm\mu$ stands for $(\lambda_1\pm \mu_1,\lambda_2\pm \mu_2,\dots)$,
and $(\lambda,\mu)$ stands for the concatenation of $\lambda$ and $\mu$. 
Note that if $\lambda \leq \mu$ and $\nu \leq \omega$, then
$\lambda \cup \nu \leq \mu \cup \omega$.  
We shall denote by $\delta_n$ (or simply $\delta$ when the value of $n$
is clear) the partition $(n-1,n-2,\dots,0)$.

Any partition $\lambda$ has an associated Ferrers diagram 
with $\lambda_i$ lattice squares in the $i^{th}$ row, 
from the bottom to top.  For example,
\begin{equation}
\lambda\,=\,(4,2)\,=\,
{\tiny{\tableau*[scY]{ & \cr & & & }}} \, .
\end{equation}
For each cell $s=(i,j)$ in the diagram of $\lambda$, let
$\ell'(s), \ell(s), a(s)$ and $a'(s)$ be respectively the number of
cells in the diagram of $\lambda$ to the south, north, east and west
of the cell $s$.  The hook-length of any cell in $\lambda$, is
$h_s(\lambda)=\ell(s)+a(s)+1$.  In the example,
$h_{(1,2)}(4,2)=2+1+1$.  The {\it main hook-length} of $\lambda$, $h_M(\lambda)$,
is the hook-length of the cell $s=(1,1)$ in the diagram of $\lambda$.  
Therefore, $h_{M}\bigl((4,2)\bigr)=5$.
The conjugate $\lambda'$ of a  partition $\lambda$ is defined
by the reflection of the Ferrers diagram about the main diagonal.
For example, the conjugate of (4,2) is
\begin{equation}
\lambda' \,=\,
{\tiny{\tableau*[scY]{ \cr \cr & \cr &   }}}
\,=\,(2,2,1,1)\,.
\end{equation}

A partition $\lambda$ is said to be $k$-{\it bounded} if its first part is 
not larger than $k$, i.e,  if $\lambda_1 \leq k$. 
We associate to any $k$-bounded partition $\lambda$
a sequence of partitions, $\lambda^{\to k}$,
called the $k$-{\it split} of $\lambda$.
$\lambda^{\to k}= (\lambda^{(1)},\lambda^{(2)},\ldots,\lambda^{(r)})$
is obtained by partitioning $\lambda$ (without rearranging
the entries) into partitions $\lambda^{(i)}$
where $h_M(\lambda^{(i)})= k$, for all $i<r $.  For example,
$(3,2,2,2,1,1)^{\to 3}=\bigl((3),(2,2),(2,1),(1)\bigr)$ 
and $(3,2,2,2,1,1)^{\to 4}=\bigl((3,2),(2,2,1),(1)\bigr)$ .
Equivalently, the diagram of $\lambda$ is cut
horizontally into partitions with main hook-length $k$.
\begin{equation}
{\tiny{\tableau*[scY]{ \cr \cr & \cr & \cr & \cr & & \cr}}} \quad
\begin{matrix}
& {\tiny{\tableau*[sbY]{  \cr}}} \hfill \\
& {\tiny{\tableau*[sbY]{  \cr & \cr}}} \hfill \\
\longrightarrow^{(3)} & {\tiny{\tableau*[sbY]{ & \cr  & \cr}}} \hfill \\
& {\tiny{\tableau*[sbY]{  & & \cr}}} \hfill
\end{matrix}
\qquad
\text{and}
\qquad
{\tiny{\tableau*[scY]{ \cr \cr & \cr & \cr & \cr & & \cr}}} \quad
\begin{matrix}
& {\tiny{\tableau*[sbY]{  \cr}}} \hfill \\
\longrightarrow^{(4)} & {\tiny{\tableau*[sbY]{\cr  & \cr  & \cr}}} \hfill \\
& {\tiny{\tableau*[sbY]{ & \cr & & \cr}}} \hfill
\end{matrix}
\, .
\end{equation}
The last partition in the sequence $\lambda^{\to k}$ may have main hook-
length less than $k$.
It is important to note that $\lambda^{\to k}=(\lambda)$
when $h_M(\lambda)\leq k$.

\subsection{Symmetric functions} The power sum $p_i(x_1,x_2,\ldots)$ is
\begin{equation}
p_i(x_1,x_2,\ldots) = x_1^i+x_2^i+\cdots \, ,
\end{equation}
and for a partition $\lambda=(\lambda_1,\lambda_2,\dots)$, 
\begin{equation}
p_{\lambda}(x_1,x_2,\ldots)=
p_{\lambda_1}(x_1,x_2,\ldots)\, p_{\lambda_2}(x_1,x_2,\ldots)\cdots \,.
\end{equation}
We employ the notation of $\lambda$-rings, needing only the formal 
ring of symmetric functions $\Lambda$ to act on the ring of rational 
functions in $x_1,\dots,x_N,q,t$, with coefficients in $\mathbb R$.
The action of a power sum $p_i$ on a rational function is, by definition,
\begin{equation}
p_{i} \left[ \frac{\sum_{\alpha} c_{\alpha} u_{\alpha} }
{ \sum_{\beta} d_{\beta} v_{\beta} } \right]
 =\frac{\sum_{\alpha} c_{\alpha} u_{\alpha}^i }
{ \sum_{\beta} d_{\beta} v_{\beta}^i},
\label{actp}
\end{equation}
with $c_{\alpha},d_{\beta} \in \mathbb R$ and $u_{\alpha}, v_{\beta}$
monomials in $x_1,\dots,x_N,q,t$.  Since the power sums form a basis
of the ring $\Lambda$, any symmetric function has a 
unique expression in terms of power sums, and \eqref{actp} extends to an 
action of $\Lambda$ on rational functions.
In particular $f[X]$, the action of a symmetric function
$f$ on the monomial $X=x_1+\cdots+x_N$, is simply $f(x_1,\ldots,x_N)$.
In the remainder of the article, we will always consider the number of
variables $N$ to be infinite, unless otherwise specified.

The monomial symmetric function $m_{\lambda}[X_n]$ is
\begin{equation}
m_{\lambda}[X_n]= \sum_{\sigma \in S_n; \sigma(\lambda) \text{~distinct}} 
x^{\sigma(\lambda)} \, .
\end{equation}
The complete symmetric function $h_r[X]$ is
\begin{equation}
h_r[X] = \sum_{1\leq i_1\leq i_2\leq\cdots\leq i_r}
x_{i_1} x_{i_2}\cdots x_{i_r} \, ,
\end{equation}
and $h_{\lambda}[X]$ stands for the homogeneous symmetric function
\begin{equation}
h_{\lambda}[X]=
h_{\lambda_1}[X]\, h_{\lambda_2}[X]\cdots \,.
\end{equation} 
And, the elementary symmetric function  $e_{\lambda}[X]$ is
defined
\begin{equation}
e_{\lambda}[X]=
e_{\lambda_1}[X]\, e_{\lambda_2}[X]\cdots \,,
\quad\text{where}\quad
e_r[X] = \sum_{1 \leq i_1< i_2<\cdots< i_r}
x_{i_1} x_{i_2}\cdots x_{i_r} \, .
\end{equation}

Although the Schur functions may be characterized in many ways, 
here it will be convenient to use the Jacobi-Trudi 
determinantal expression:
\begin{equation}
s_\lambda[X] = \text{det}\Bigl | h_{\lambda_i+j-1}[X] 
\Bigr|_{1 \leq i,j \leq \ell(\lambda)}
\label{jt}
\end{equation}
where $h_r[X]=0$ if $r < 0$.  
Note, in particular, $s_r[X]=h_r[X]$.

We recall that the Macdonald scalar product,  
$\langle \ , \ \rangle_{q,t}$, 
on $\Lambda \otimes \mathbb Q(q,t)$ is defined by setting 
\begin{equation}
\langle p_\lambda[X], p_\mu[X] \rangle_{q,t}
        =\delta_{\lambda \mu } \, z_\lambda \prod_{i=1}^{\ell(\lambda)} \frac{1-q^{\lambda_i}}{1-t^{\lambda_i}} \, ,
\end{equation}
where for a partition $\lambda$ with $m_i(\lambda)$ parts equal to $i$,
we associate the number
\begin{equation} \label{1}
z_\lambda
        = 1^{m_1} m_1 !  \, 2^{m_2} m_2! \dotsm
\end{equation}
If $q=t$, this expression no longer depends
on a parameter and is then denoted $\langle, \rangle$,  
satisfying
\begin{equation} \label{schurscalar}
\langle s_{\lambda}[X],s_{\mu}[X]\rangle= \delta_{\lambda \mu} \, .
\end{equation} 
The Macdonald integral forms 
$J_\lambda [X; q,t]$ are uniquely characterized \cite{[Ma]} by
\begin{align}
  \mathrm{(i)} \ &  \langle J_\lambda, J_\mu \rangle_{q,t} = 0, \qquad \text{if } \lambda \ne \mu , \\
  \mathrm{(ii)} \ &  J_\lambda[X;q,t] = \sum_{\mu \le \lambda} v_{\lambda\mu }(q,t) 
s_\mu[X]  
\quad\text{with}\quad v_{\lambda \mu}(q,t) \in \mathbb Q(q,t)\,,  \\
  \mathrm{(iii)} \ & v_{\lambda\lambda}(q,t)= 
\prod_{s \in \lambda} (1-q^{a(s)} t^{\ell (s)+1}),
\end{align}
Here, we use a modification of the Macdonald 
integral forms that is obtained by setting
\begin{equation}
H_{\lambda}[X;q,t]=J_{\lambda}[X/(1-t);q,t]= \sum_{\mu} K_{\mu \lambda}(q,t)\, 
s_{\mu}[X] \, ,
\end{equation}
with the coefficients $K_{\mu \lambda}(q,t) \in \mathbb N[q,t]$ 
known as the $q,t$-Kostka polynomials.
When $q=0$, $J_{\lambda}[X;q,t]$ reduces to the Hall-Littlewood
polynomial, $J_{\lambda}[X;0,t]=Q_{\lambda}[X;t]$.  Again, we 
use a modification;
\begin{equation}
H_{\lambda}[X;t]=H_{\lambda}[X;0,t]=Q_{\lambda}[X/(1-t);t]= 
s_{\lambda}[X]+\sum_{\mu> \lambda} K_{\mu \lambda}(t)\, 
s_{\mu}[X] \, ,
\label{HallinS}
\end{equation}
with the coefficients $K_{\mu \lambda}(t)\in\mathbb N[t]$ known 
as the Kostka-Foulkes polynomials.
The Kostka numbers $K_{\mu \lambda} \in \mathbb N$ arise
in the limit $t=1$, as coefficients in the expansion
\begin{equation} \label{kostka}
h_{\lambda}[X]=H_{\lambda}[X;1]= s_{\lambda}[X]+\sum_{\mu>\lambda} 
K_{\mu \lambda}\, s_{\mu}[X] \, .
\end{equation}
The Kostka numbers also appear in the expansion
\begin{equation} \label{kostkam}
s_{\lambda}[X] = m_{\lambda}[X]+\sum_{\mu<\lambda} 
K_{\lambda \mu}\, m_{\mu}[X] \, .
\end{equation}

\section{Vertex operators and Schur functions}
\label{SecVer}
The ring of symmetric polynomials over rational functions in  
the parameter $t$ has shown to be of interest in many fields 
of mathematics and physics.  A natural basis for this space 
is given by the Hall-Littlewood polynomials, $H_{\lambda}[X;t]$; 
$t$-analogs of the homogeneous symmetric functions, $h_{\lambda}[X]$.  
Our approach to the study of this space employs vertex operators 
that arise in the recursive construction for the Hall-Littlewood 
polynomials \cite{[Ji]}.  These operators can be 
defined \cite{[SZ]} for $\ell \in \mathbb Z$, by 
\begin{equation} 
\label{vertexop}
B_\ell = \sum_{i=0}^{\infty} s_{i+\ell}[X] \, s_{i}[X(t-1)]^{\perp} \, , 
\end{equation}
where for $f,g$ and $h$ arbitrary symmetric functions,  
$f^{\perp}$ is such that on the scalar product \eqref{schurscalar},
\begin{equation}
\langle f^{\perp} g,h \rangle = \langle g,fh \rangle \, .
\end{equation}
Note, $B_{\ell} \cdot 1=0$ if $\ell<0$.
The operators build the Hall-Littlewood polynomials by
\begin{equation} 
H_\lambda[X;t] = B_{\lambda_1} H_{\lambda_2,\ldots,\lambda_\ell}[X;t]
\, , \qquad \text{for} \qquad \lambda_1 \geq \lambda_2 \, ,
\end{equation} 
and, since they satisfy the relation
\begin{equation}
\label{commuB}
B_m B_n = t B_n B_m +t B_{m+1} B_{n-1}- B_{n-1} B_{m+1} \, ,
\qquad m,n \in \mathbb Z \, , 
\end{equation}
their action on 
$\Lambda$ can be computed algebraically.

The definition for these operators 
was extended in \cite{[SZ]} to any partition 
$\lambda$ of length $L$, by
\begin{equation}
\label{jtt}
B_\lambda \equiv \prod_{1\leq i<j\leq L}
(1-te_{ij})
B_{\lambda_1} \cdots
B_{\lambda_L} \, ,
\end{equation}
where $e_{ij}$ acts by
\begin{equation}
e_{ij}\left(
B_{\lambda_1} \cdots
 B_{\lambda_L}
\right)
=
B_{\lambda_1} \cdots
B_{\lambda_i+1}
\cdots
B_{\lambda_j-1} \cdots
B_{\lambda_L}
\, .
\end{equation}

It is important to note that \eqref{vertexop} gives 
$B_\ell=s_\ell$ when $t=1$, and thus \eqref{jtt} 
reduces to the Jacobi-Trudi formula \eqref{jt} for $s_\lambda[X]$.  Therefore,
we have
\begin{equation} \label{Bt1}
B_{\lambda}=s_{\lambda}[X] \, , \qquad \text{when~} t=1 \, .
\end{equation}

\subsection{Identities on Schur functions and vertex operators}
Here we derive properties for the vertex operators,
and consequently for the Schur functions.  
Our properties reveal information about the behavior of 
$B_\lambda$ and as a by-product, allow us to prove conjectures
relating to a filtration of the symmetric function
space.
The main result of this section is an explicit formula for
special products of two vertex operators. 

\bigskip

\begin{theorem} \label{theorembig}
Let $a,r$ and $m$ be nonnegative integers
and $\nu$ be a partition with $\ell(\nu)\leq r$. 
Then
\begin{equation}
B_{a^{r+m}} \,  B_{\nu}
 = \sum_{\mu\in\P_{\leq m}^{r}} 
(-t)^{|\mu|} B_{a^r+\mu} \, B_{a^m-(\mu')^R,\nu}\, ,
\end{equation}
where the summand vanishes if $(a^m-(\mu')^R)$
is not a partition.
\end{theorem}

\bigskip

Since ${B_\lambda}\to s_\lambda$ when $t=1$ by (\ref{Bt1}),
an immediate consequence of this result is
an explicit formula for a product of  two Schur functions.

\begin{corollary} 
Let $a,r$ and $m$ be nonnegative integers
and $\nu$ be a partition with $\ell(\nu)\leq r$. 
Then
\begin{equation}
s_{a^{r+m}} \,  s_{\nu}
 = \sum_{\mu\in\P_{\leq m}^{r}} 
(-1)^{|\mu|} s_{a^r+\mu} \, s_{a^m-(\mu')^R,\nu}\, ,
\end{equation}
where the summand vanishes if $(a^m-(\mu')^R)$
is not a partition.
\end{corollary}

From Theorem~\ref{theorembig}, we derive a number of 
identities for the product of vertex operators.
We now state and prove these identities, postponing the 
proof of our theorem to the end of this section.  
We start by giving several properties of the vertex operators $B_v$
that arise from the fact that they satisfy the same reordering
relations as the Schur functions.  Namely
(see Proposition~3 of \cite{[SZ]}),
for an integral vector $v$, 
\begin{equation}
B_{v} = -B_{v_1,\ldots,v_{i-1},v_{i+1}-1,v_i+1,v_{i+2},\ldots,v_r}
\,.
\label{commute}
\end{equation}
This reordering relation allows us to prove the following properties:

\bigskip

\begin{property} \label{lemzero}
If $v=(v_1,\dots,v_r)\in\mathbb Z^r$ is such that  
$v_j-v_i=j-i$ for $i\ne j$,  then $B_v=0$.
\end{property}
\noindent{\bf Proof.}  
Without loss of generality, assume $j>i$.
Given $B_v$, we successively apply \eqref{commute}
to move $v_j$ from position $j$ to position $i+1$;
\begin{eqnarray}
B_v & = & (-1)^{j-i-1}
B_{v_1,\dots,v_i,v_{j}-(j-i-1),v_{i+1}+1,\dots,v_{j-1}+1,v_{j+1},\dots,v_r}
\\
& = &(-1)^{j-i-1}
B_{v_1,\dots,v_i,v_{i}+1,v_{i+1}+1,\dots,v_{j-1}+1,v_{j+1},\dots,v_r} \, ,
\end{eqnarray}
since $v_j-v_i=j-i$. 
Further, switching the entries in position $i$ and $i+1$, we obtain
\begin{equation}
 B_{v_1,\dots,v_i,v_{i}+1,v_{i+1}+1,\dots,v_{j-1}+1,v_{j+1},\dots,v_r} 
= -
B_{v_1,\dots,v_i,v_{i}+1,v_{i+1}+1,\dots,v_{j-1}+1,v_{j+1},\dots,v_r} \, ,
\end{equation}
which implies that 
$B_{v_1,\dots,v_i,v_{i}+1,v_{i+1}+1,\dots,v_{j-1}+1,v_{j+1},\dots,v_r}=0$.
Consequently, $B_v$ is also null. 
\hfill $\square$

\bigskip

\begin{property} \label{lemmax}
Let $\mu$ and $\nu$ be partitions of lengths $m$ and $r$ respectively.
If there is a non-zero partition reordering, $\pm B_\lambda$,
of $B_{\mu,\nu}$, then $\lambda_1 = \max \{\mu_1,\nu_1-m \}$ and 
$\lambda_{m+r} = \min \{\mu_m+r,\nu_r \}$. 
\end{property}
\noindent{\bf Proof.}\quad 
An element of $v$ in $B_v$ moved $i$ steps to the left (right)
with (\ref{commute}), is decreased (increased) by $i$.  
Therefore, the only possible entries in the first position of
any reordering of $B_{\mu,\nu}$ are
\begin{equation}
\mu_i-(i-1) \, , \quad i=1,\dots,m  \quad\text{or} \quad 
\nu_j-(j+m-1) \, , \quad j=1,\dots,r \, .
\end{equation}
Assume the first entry in a reordering $B_\lambda$ 
of $B_{\mu,\nu}$ is not the largest of these,
i.e.  $\lambda_1<\max\{\mu_1,\nu_1-m \}$.
Since the entries $\mu_1+i$ for some $i\geq 0$ and $\nu_1-j$ 
for some $j\leq m$ must occur in $B_{\lambda}$,
and $\lambda_1<\max\{\mu_1,\nu_1-m \}\leq\max\{\mu_1+i,\nu_1-j\} $,
$\lambda_1$ cannot be the largest entry.
Therefore $\lambda$ is not a partition unless 
$\lambda_1=\max \{\mu_1,\nu_1-m \}$.
Similar reasoning applies for the smallest entry $\lambda_{m+r}$. 
\hfill $\square$

The contrapositive of Property~\ref{lemmax} then 
gives the following result:

\bigskip

\begin{property} \label{lemmax2}
Let $\mu$ and $\nu$ be partitions of lengths $m$ and $r$, 
respectively.  If $\max \{\mu_1,\nu_1-m \} < \min \{\mu_m+r,\nu_r \}$
then $B_{\mu,\nu}=0$.
\end{property}

These properties of $B_v$ allow us to give several identities 
concerning the product of an operator $B_\nu$ with
an operator indexed by partitions of the form $(\ell^{k+1-\ell})$,
$\ell=1,\dots,k$, 
hereafter referred to as $k$-rectangles.  
Our identities are derived from
particular cases of  Theorem~\ref{theorembig}.  
The first is a $t$-commutation relation.

\bigskip

\begin{identity} \label{identityzero}
Let $i$, $k$ and $\ell$ be nonnegative integers.  For 
$\ell \leq i \leq k$, we have 
\begin{equation}
B_{\ell^{k+1-\ell}} \, B_{i} = t^{i-\ell} \, B_{i} \, B_{\ell^{k+1-\ell}} \, .
\end{equation}
\end{identity}
\noindent {\bf Proof.} \quad 
Theorem~\ref{theorembig}, 
with $a=\ell$, $r=1$, $m=k-\ell$, and $\nu=(i)$, gives
\begin{equation} \label{eqtcommute}
\begin{split}
B_{\ell^{k+1-\ell}} \, B_{i}& = 
\sum_{\mu\in\P^1_{\leq k-\ell} }
(-t)^{|\mu|}
B_{\ell+\mu}  \, B_{\ell^{k-\ell}-(\mu')^R,i} 
 = \sum_{j=0}^{k-\ell} (-t)^j B_{\ell+j} \, B_{\ell^{k-\ell}-(1^j)^R,i} 
\, .
\end{split}
\end{equation}
Since $(\ell^{k-\ell}-(1^j)^R,i) =(\ell^{k-\ell-j},(\ell-1)^j,i)$,
it suffices to show $B_{\ell^{k-\ell-j},(\ell-1)^j,i} =0$ 
unless $j=i-\ell$, in which case
$B_{\ell^{k-i},(\ell-1)^{i-\ell},i}=(-1)^{i-\ell}B_{\ell^{k+1-\ell}}$
by relation \eqref{commute}. 
When $k \geq i>j+\ell$, positions $k-\ell+1$ and $k-i+1$ 
of $(\ell^{k-\ell-j},(\ell-1)^j,i)$ contain the entries
$i$ and $\ell$, resp. Since these entries differ by $i-\ell$, 
Property~\ref{lemzero} gives  $B_{\ell^{k-\ell-j},(\ell-1)^j,i}=0$.
Similarly, when $\ell \leq i<\ell+j$, 
the entries in positions $k-\ell+1$ and $k-i$, resp
$i$ and $\ell-1$, differ by $i-\ell+1$
and again we have $B_{\ell^{k-\ell-j},(\ell-1)^j,i}=0$.
\hfill $\square$

\bigskip

\begin{identity} \label{identity4}
Let $\nu$ be a partition of length $r$ where $h_M(\nu)\leq k$
and $\nu_1\geq\ell$ for $\ell\in \mathbb N$.  Then
\begin{eqnarray}
B_{\ell^{k+1-\ell}} \, B_{\nu}
=
\sum_{i} (-t)^{c(i)} B_{\rho^{(i)}} \, B_{\gamma^{(i)}}
\, ,
\end{eqnarray}
where $c(i)\in\mathbb N$, and 
$\rho^{(i)}$ and $\gamma^{(i)}$ are partitions
such that 
$\rho^{(i)}_1=\nu_1$, $\rho^{(i)}_L \geq \ell$, 
$h_M(\rho^{(i)})=k, \gamma^{(i)}_1=\ell$,
and $h_M(\gamma^{(i)})\leq k$ (equality holds only
when $h_M(\nu)= k$).
\end{identity}
\noindent {\bf Proof.} \quad 
Let $m=\nu_1-\ell$, $r=k+1-\nu_1$, and $a=\ell$ in 
Theorem~\ref{theorembig}.  We then have
\begin{equation}
\label{eqzero}
B_{\ell^{k+1-\ell}} \, B_{\nu}
 = \sum_{\ell(\mu)\leq k+1-\nu_1\atop
\mu_1\leq\nu_1-\ell}
(-t)^{|\mu|} B_{\ell^{k+1-\nu_1}+\mu} \, B_{\ell^{\nu_1-\ell}-(\mu')^R,\nu}\, .
\end{equation} 
If $\mu_1<\nu_1-\ell$ (equivalently $\mu'_{\nu_1-\ell}=0$),
then the first and $\nu_1-\ell+1^{th}$ entries
of $({\ell^{\nu_1-\ell}-(\mu')^R,\nu})$
are $\ell-\mu_{\nu_1-\ell}'=\ell$ and $\nu_1$ respectively.
Since these entries differ by $\nu_1-\ell$,
Property~\ref{lemzero} implies that
$B_{\ell^{\nu_1-\ell}-(\mu')^R,\nu}=0$ in this case.
We thus have, 
\begin{equation}
\label{te1}
B_{\ell^{k+1-\ell}} \, B_{\nu}
 =\sum_{ \ell(\mu)\leq k+1-\nu_1 \atop \mu_1=\nu_1-\ell }
(-t)^{|\mu|} B_{\ell^{k+1-\nu_1}+\mu} \, B_{\ell^{\nu_1-\ell}-(\mu')^R,\nu}
\, .
\end{equation} 
If we let $(\ell^{k+1-\nu_1}+\mu)=\rho^{(i)}$
then $\rho_1^{(i)}=\mu_1+\ell=\nu_1$.
Moreover,  $\ell(\mu)\leq k+1-\nu_1$
gives that $\rho_{L}^{(i)}\geq\ell$ and $h_M(\rho^{(i)})=\nu_1+(k+1-\nu_1)-1=k$.
Further, since $(\ell^{\nu_1-\ell}-(\mu')^R,\nu)$ is
the concatenation of two weakly decreasing sequences,
Property~\ref{lemmax} implies any non-zero 
partition reordering, $\pm B_{\gamma^{(i)}}$,
of $B_{\ell^{\nu_1-\ell}-(\mu')^R,\nu}$ 
has $\gamma^{(i)}_1=\max \{\ell-\mu'_{\nu_1-\ell}, \nu_1-(\nu_1-\ell) \}=\ell$. 
Then, $h_M(\gamma^{(i)})=\ell+(\nu_1-\ell+r)-1=\nu_1+r-1\leq k$ 
since $h_M(\nu)=\nu_1+r-1\leq k$ (equality holds only when
$h_M(\nu)=k$).
\hfill $\square$

\medskip

We now give two identities regarding the product 
of an operator indexed by a $k$-rectangle, $(\ell^{k+1-\ell})$,
with an operator indexed by any partition with main 
hook-length exactly $k$.

\bigskip

\begin{identity}  \label{identity2}
Let $\nu$ be a partition of length $r$ where $h_M(\nu)=k$,
$\nu_1\geq\ell$, and $\nu_r < \ell$.  Then
\begin{eqnarray}
B_{\ell^{k+1-\ell}} \, B_{\nu}
=
\sum_{i} (-t)^{c(i)} B_{\rho^{(i)}} \, B_{\gamma^{(i)}} \, ,
\end{eqnarray}
where $c(i)\in\mathbb N$, and 
$\rho^{(i)}$ and $\gamma^{(i)}$ are partitions
where
$\rho^{(i)}_1=\nu_1$, $\rho^{(i)}_L \geq \ell$, 
$h_M(\rho^{(i)})=h_M(\gamma^{(i)})=k$, 
$\gamma^{(i)}_1=\ell$, and $\gamma^{(i)}_L=\nu_r$. 
\end{identity}
\noindent{\bf Proof.}\quad This result follows from 
formula (\ref{te1}) in Identity~\ref{identity4} 
with $h_M(\nu)=\nu_1+r-1=k$ and $\nu_r<\ell$. 
Since $\ell(\mu)\leq k+1-\nu_1=r$ implies $\mu_1'\leq r$,
Property~\ref{lemmax} proves 
$\gamma_{L}^{(i)}=\min \{\ell-\mu_1'+r, \nu_r \}=\nu_r$ 
since $\ell-\mu_1'+r\geq \ell-r+r=\ell>\nu_r$.
\hfill $\square$

\begin{lemma} \label{lemunique}
Let $\mu$, $\nu$, and $\lambda$ be partitions, and  let 
$\lambda = w\bigl( (\mu,\nu)+\delta) \bigr)-\delta$ for some permutation
$w$.  Let $\gamma$ be a partition such that $\gamma\not = \mu$, and 
such that $\ell(\mu)=\ell(\gamma)=n$, where both partitions may contain zeroes.
Then $\lambda \neq \sigma\bigl( (\gamma,\nu)+\delta) \bigr)-\delta$
for all permutations $\sigma$.
\end{lemma}
\noindent {\bf Proof.} \quad 
Assume there exist permutations 
$w$ and $\sigma$ where 
$\lambda= w\bigl( (\mu,\nu)+\delta) \bigr)-\delta=
\sigma\bigl( (\gamma,\nu)+\delta) \bigr)-\delta$
for $\mu\neq\gamma$. 
This implies that
$(\mu,\nu)+\delta= \sigma'\bigl( (\gamma,\nu)+\delta) \bigr)$ for some
permutation $\sigma'$.
For all $i>n$, 
$\bigl((\mu,\nu)+\delta\bigr)_i=\bigl((\gamma,\nu)+\delta\bigr)_i$,
since $\ell(\gamma)=\ell(\mu)=n$ implies $\gamma_i=\mu_i=0$. 
Thus, $\mu+\delta_n=\sigma'(\gamma+\delta_n)$, 
or equivalently $\mu=\sigma'(\gamma+\delta_n)-\delta_n$, 
for some $\sigma'\in S_n$.
However, since $\gamma$ is the only element that is a partition in the set 
$\{\tau(\gamma+\delta_n)-\delta_n\,:\, \tau\in S_n\}$, 
we arrive at the contradiction $\mu=\gamma$. 
\hfill $\square$

\bigskip

\begin{identity}  \label{identity3}
Let $\nu$ be a partition of length $r$ where $h_M(\nu)=k$
and $\nu_r \geq \ell$.  Then
\begin{equation}
B_{\ell^{k+1-\ell}} \, B_{\nu}
= t^{|\nu|-r \ell } B_{\nu} \, B_{\ell^{k+1-\ell}}\,.
\end{equation} 
\end{identity}
\noindent{\bf Proof.}\quad 
For $\gamma$ indexing the partition reordering of
$B_{\ell^{\nu_1-\ell}-(\mu')^R,\nu}$,
Identity~\ref{identity4} implies
that \eqref{te1} holds with $\gamma_1=\ell$.
Since $h_M(\nu)=\nu_1+r-1=k$, we have
$\ell(\gamma^{(i)})=\nu_1-\ell+r=k+1-\ell$.  Also, 
$\ell(\mu) \leq k+1-\nu_1=r$ leads to  $\mu_1' \leq r$.
Therefore, $\nu_r\geq\ell$ implies by Property~\ref{lemmax} 
that $\gamma_{k+1-\ell}=\min \{\ell-\mu_1'+r, \nu_r \}\geq\ell$ 
since $\ell-\mu_1'+r\geq \ell-r+r=\ell$.
Thus, from Property~\ref{lemmax2},
$\gamma_1=\gamma_{k+1-\ell}=\ell$ implies 
that all non-zero $B_{\ell^{\nu_1-\ell}-(\mu')^R,\nu}$
equal $\pm B_{\gamma}$ where $\gamma=(\ell^{k+1-\ell})$.

It now suffices to show that 
$B_{\ell^{\nu_1-\ell}-(\mu')^R,\nu} 
=\pm B_{(\ell^{k+1-\ell})}$
only when $\mu=(\nu_1-\ell,\dots,\nu_r-\ell)$. 
Our claim will follow since
\eqref{te1} then simplifies to 
\begin{equation}
\label{signe}
B_{\ell^{k+1-\ell}} \, B_{\nu}
 = \pm \, t^{|\nu|-r \ell } B_{\nu} \, B_{\ell^{k+1-\ell}}\, ,
\end{equation} 
where the sign must be positive since 
when $t=1$, this relation becomes
$s_{\ell^{k+1-\ell}} s_{\nu}=s_{\nu}s_{\ell^{k+1-\ell}}$.
In fact, we only need to show that
when $\mu=(\nu_1-\ell,\dots,\nu_r-\ell)$,
there exists some permutation $w$ where
\begin{equation}
\label{eqlem}
w\big((\ell^{\nu_1-\ell}-(\mu')^R,\nu)+\delta\big)-\delta
=(\ell^{k+1-\ell})\,.
\end{equation}
Then by Lemma~\ref{lemunique},
since there exists no other partition 
$(\ell^{\nu_1-\ell}-(\mu')^R)$ such that 
$w\big((\ell^{\nu_1-\ell}-(\mu')^R,\nu)+\delta\big)-\delta
=(\ell^{k+1-\ell})$, non-zero terms occur only when
$\mu=(\nu_1-\ell,\dots,\nu_r-\ell)$.  

To prove \eqref{eqlem}, it is equivalent to 
show that there exists some permutation $w'$ where
\begin{equation}
w'\delta+(\ell^{k+1-\ell})
=
(\ell^{\nu_1-\ell}-(\mu')^R,\nu)+\delta\big)
\end{equation}
Since $\mu=(\nu_1-\ell,\dots,\nu_r-\ell)
\implies \nu=(\mu_1+\ell,\dots,\mu_r+\ell)$,
and $\delta=(k-\ell,\dots,0)=(r+\mu_1-1,\dots,0)$
given $r=k-\nu_1+1$ (from $h_M(\nu)=k$),
we must show there is some permutation $w'$ where
\begin{equation}
\label{eq2lem}
\begin{split}
 w'\delta+(\ell^{k+1-\ell})
& 
= 
(\ell-\mu_{\mu_1}'+r+\mu_1-1,\dots,\ell-\mu_1'+r,\mu_1+\ell+r-1,
\dots,\mu_r+\ell) 
\\
 & 
= 
(\ell^{k+1-\ell})+(-\mu_{\mu_1}'+r+\mu_1-1,
\dots,-\mu_1'+r,\mu_1+r-1,\dots,\mu_r)
 \, .
\end{split}
\end{equation}
The last $r$ entries of 
$u=(-\mu_{\mu_1}'+r+\mu_1-1,\dots,-\mu_1'+r,\mu_1+r-1,\dots,\mu_r)$
are $\mu_i+r-i$, $i=1,\dots,r$, and the first $\mu_1$ entries
are $r-1+j-\mu_j'$, $j=\mu_1,\dots,1$.  Since $r \geq \mu'_1$,
a vector of this type is known \cite{[Ma]} to be a permutation 
of $\delta_{r+\mu_1}$.  Thus, there is some 
$w'$ such that $u=w'\delta$, and \eqref{eq2lem} follows. \hfill $\square$

\subsection{Proof of Theorem~\ref{theorembig}}
We use several lemmas that rely on 
properties for the Kostka matrix. These properties
are derived in Appendix~\ref{A1} for lack of reference.
Here, we use $E_m^{d}$ to denote the set of vectors of length $m$, 
with $d$ ones and $m-d$ zeroes. Then
$E_m^{\lambda}$ is the set of vectors $v=v_1+v_2+\dots$, where
$v_1 \in E_m^{\lambda_1}$, $v_2 \in E_m^{\lambda_2},\ldots$.

\begin{lemma} \label{iden1}
For a partition $\lambda \in \P^r$ and any partition $\nu$, we have
\begin{equation}
\sum_{\sigma \in S_r;\sigma(\lambda)~{\rm{ distinct}}}
B_{\sigma(\lambda)+b^r,\nu} = \sum_{\mu \in {\P}^r} 
B_{\mu+b^r,\nu} K_{\lambda \mu}^{-1} \, .
\end{equation}
\end{lemma}
\noindent{\bf{Proof.}}\quad 
Using Formula~\ref{monoschur}, we have
\begin{equation}
\label{mac}
m_{\lambda+b^r}= \sum_{\sigma \in S_r; \sigma(\lambda+b^r)~{\rm{ distinct}}} 
s_{\sigma(\lambda+b^r)}
= \sum_{\mu \in {\mathcal P}^r} K_{\lambda+b^r\,
\mu+b^r}^{-1} \,s_{\mu+b^r} \, .
\end{equation}
With
$K_{\lambda\mu}^{-1}=K^{-1}_{\lambda+b^r\mu+b^r}$, from Formula~\ref{kostkaplus},
 this gives 
$$
\sum_{\sigma \in S_r; \sigma(\lambda)~{\rm{ distinct}}} 
s_{\sigma(\lambda)+b^r}
= \sum_{\mu\in\mathcal P^r} K_{\lambda\mu}^{-1} \,s_{\mu+b^r}\, .
$$
We formally replace $s_*$ by $B_{*,\nu}$ since $s_{\sigma(\lambda)}$ and 
$B_{\sigma(\lambda),\nu}$ both obey the reordering relation 
\eqref{commute}.  \hfill$\square$

\begin{lemma} \label{iden2}
For $\lambda \in \P^{m}_{\leq a}$,
we have
\begin{equation}
\sum_{E \in E^{\lambda}_m} B_{a^m-E,\nu} =
\sum_{\omega}\sum_{\gamma\in \P^m} \sum_{\rho \in P^m_{\leq a}}  K_{\gamma \rho}^{-1}
K_{\omega \lambda} K_{\omega' \gamma} \, B_{a^m-\rho^R,\nu} \, .
\end{equation}
\end{lemma}
\noindent{\bf{Proof.}}
It is known \cite{[Ma]} that for $\lambda \in \P^m$, 
\begin{equation}
e_{\lambda}= \sum_{\gamma,\omega} K_{\omega \lambda} K_{\omega' \gamma} 
m_{\gamma} \quad\text{and}\quad
\sum_{E \in E^{\lambda}_m} x^E =
e_{\lambda}[x_1+\cdots+x_m] \,.
\end{equation}
Equivalently, we thus have
\begin{equation}
\sum_{E \in E^{\lambda}_m} x^E = 
\sum_{\omega} \sum_{\gamma \in \P^m} K_{\omega \lambda} K_{\omega' \gamma} \, 
\sum_{\sigma \in S_m;\sigma(\gamma)~{\rm{distinct}}} x^{\sigma(\gamma)} \, .
\end{equation}
Formally replacing $x^*$ by $B_{a^m-*,\nu}$, this implies
\begin{equation}
\sum_{E \in E^{\lambda}_m} B_{a^m-E,\nu}= 
\sum_{\omega} \sum_{\gamma \in \P^m} K_{\omega \lambda} K_{\omega' \gamma} \, 
\sum_{\sigma \in S_m;\sigma(\gamma)~{\rm{distinct}}} B_{a^m-\sigma(\gamma),\nu} \, .
\end{equation}
Since $a^m-\gamma^R$ is a partition, we can use Lemma~\ref{iden1}
with $\lambda=a^m-\gamma^R$ and $b=0$ to obtain
\begin{equation}
\sum_{E \in E^{\lambda}_m} B_{a^m-E,\nu} = 
\sum_{\omega} \sum_{\gamma \in \P^m} K_{\omega \lambda} K_{\omega' \gamma} \, 
\sum_{\mu \in \P^m} B_{\mu,\nu} K_{a^m-\gamma^R,\mu}^{-1}.
\end{equation}
Since $K^{-1}_{a^m-\gamma^R,\mu}=0$ if $a^m-\gamma^R \not \geq \mu$,
we have $a \geq \mu_1$. Thus, $\mu = a^m-\rho^R$ for 
some $\rho \in \P^{m}_{\leq a}$ and
\begin{equation}
\sum_{E \in E^{\lambda}_m} B_{a^m-E,\nu} = 
\sum_{\omega} \sum_{\gamma \in \P^m} K_{\omega \lambda} K_{\omega' \gamma} \, 
\sum_{\rho \in \P^m_{\leq a}} B_{a^m-\rho^R,\nu} K_{a^m-\gamma^R,a^m- \rho^R}^{-1}.
\end{equation}
The property then follows from Formula~\ref{kostkamoins}, which gives
$K_{a^m-\gamma^R,a^m- \rho^R}^{-1}=K_{\gamma,\rho}^{-1}$.
\hfill$\square$

The last lemma needed to derive our expression for the product 
of Schur functions and operators uses methods presented in \cite{[Ad]} 
and \cite{[SZ]}.  Here, $v\in[n]$ is a vector with 
entries from $0,1,\ldots,n$.
\begin{lemma}
\label{genlemma}
Let $\nu,\mu,\gamma$ be any partitions with
$\ell(\nu)=n$, $\ell(\mu)=r$ and $\ell(\gamma)=m$.
Then
\begin{equation}
\sum_{I=(i_1,\ldots,i_{m})\in[n]}
\sum_{E\in E^{I}_n}
(-t)^{|I|}
B_{\mu,\gamma+{I}} \,
B_{\nu-E}
=
\sum_{I=(i_1,\ldots,i_{r})\in [m]}
\sum_{E\in E^{I}_{m}}
(-t)^{|I|}
B_{\mu+I} \,
B_{\gamma-E,\nu} \, .
\end{equation}
\end{lemma}
\noindent{\bf Proof.}
Identity (20) in \cite{[SZ]} implies that
\begin{equation}
\label{oneeq}
H(U^r,V^m)H(Z^n)
\prod_{i=1}^m
\prod_{j=1}^n
\left( 1-\frac{t}{v_i} z_j\right)
=
H(U^r)H(V^m,Z^n)
\prod_{i=1}^r
\prod_{j=1}^m
\left( 1-\frac{t}{u_i} v_j\right)
\, ,
\end{equation}
where $H(X^\ell)$ is a formal Laurent series in an
ordered set of variables $X^\ell=(x_1,\ldots,x_\ell)$
with coefficients given by operators which act
on  $P\in \Lambda$ by 
\begin{equation}
H(X^\ell)P[Y] = P\left[Y-(1-q)\sum_{i=1}^\ell x_i^{-1}\right]
\prod_{i,j}\frac{1}{(1-y_ix_j)}\prod_{1\leq i<j\leq \ell}
\left(1-x_j/x_i \right)
\, .
\end{equation}
Binomial expansion of \eqref{oneeq} then gives
$$
H(U^r,V^m)H(Z^n)
\prod_{i=1}^m
\sum_{0\leq l_i \leq n}
\left( \frac{-t}{v_i} \right)^{l_i}
\sum_{E\in E_n^{l_i}}
z^E
= H(U^r)H(V^m,Z^n)
\prod_{i=1}^r
\sum_{0\leq l_i \leq m}
\left( \frac{-t}{u_i} \right)^{l_i}
\sum_{E\in E_m^{l_i}}
v^E
\, .
$$
Equivalently, by expanding the products, we have
$$
H(U^r,V^m)H(Z^n)
\sum_{0\leq l_1 \leq n}
\cdots
\sum_{0\leq l_m \leq n}
\frac{(-t)^{l_1+\cdots+l_m}}
{v_1^{l_1}\cdots v_m^{l_m}}
\sum_{E\in E_n^{(l_1,\cdots,l_m)}}
z^E
\qquad \qquad \qquad
$$
$$
\qquad \qquad \qquad
= H(U^r)H(V^m,Z^n)
\sum_{0\leq l_1 \leq m}
\cdots
\sum_{0\leq l_r \leq m}
\frac{(-t)^{l_1+\cdots+l_r}}
{u_1^{l_1}\cdots u_r^{l_r}}
\sum_{E\in E_m^{(l_1,\cdots,l_r)}}
v^E
$$
It is known \cite{[SZ]} that 
$H(X^\ell)P[Y]\big|_{x^\lambda}=B_\lambda P[Y]$.
Therefore, for partitions $\mu,\gamma,\nu$ with
$\ell(\mu)=r,\ell(\gamma)=m,\ell(\nu)=n$,
we take the coefficient of $u^\mu$, $v^\gamma$ and $z^\nu$ in
both sides to obtain
$$
\sum_{0\leq l_1 \leq n}
\cdots
\sum_{0\leq l_m \leq n}
(-t)^{l_1+\cdots+l_m}
B_{\mu,\gamma_1+l_1,\ldots,\gamma_m+l_m} 
\sum_{E\in E_n^{(l_1,\ldots,l_m)}}
B_{\nu-E}
\qquad \qquad \qquad
$$
$$
\qquad \qquad \qquad
=
\sum_{0\leq l_1 \leq m}
\cdots
\sum_{0\leq l_r \leq m}
(-t)^{l_1+\cdots+l_r}
B_{\mu_1+l_1,\cdots,\mu_r+l_r}
\sum_{E\in E_m^{(l_1,\ldots,l_r)}}
B_{\gamma-E,\nu}
\, .
$$
This completes the proof.\hfill $\square$

We can now prove Theorem~\ref{theorembig} using 
these lemmas and our properties for the vertex operators.

\noindent {\bf Proof of Theorem~\ref{theorembig}.} 
With $\ell(\nu)=n \leq r$, letting $\mu=(a^r)$ and $\gamma=(a^m)$ in
Lemma~\ref{genlemma}, we have 
\begin{equation}
\sum_{L=(l_1,\ldots,l_m)\in[n]} 
\sum_{E\in E_n^{L}}
(-t)^{|L|}
B_{a^r,a+l_1,\ldots,a+l_m} 
B_{\nu-E} 
= 
\! \! \!
\!
\sum_{L=(l_1,\ldots,l_r)\in[m]}
\sum_{E\in E_m^{L}}
(-t)^{|L|}
B_{a+l_1,\cdots,a+l_r}
B_{a^m-E,\nu}
\,.
\end{equation}
Property~\ref{lemzero} implies
$B_v=B_{a^r,a+l_1,\ldots,a+l_m}=0$ if
$1\leq l_1\leq r$, since $v_{r+1}-v_{r+1-l_1}=l_1$.
Thus $B_v=0$ unless $l_1$ is zero
since $0\leq l_1\leq n\leq r$.
Given $l_1=0$, Property~\ref{lemzero} implies
$B_v=B_{a^{r+1},a+l_2,\ldots,a+l_m}=0$ if
$1\leq l_2\leq r$, since $v_{r+2}-v_{r+2-l_2}=l_2$.
Thus $B_v=0$ unless $l_2$ is zero
since $0\leq l_2\leq n\leq r$. 
Repeating this argument, $B_{a^r,a+l_1,\ldots,a+l_m}=0$ 
unless $l_1=l_2=\cdots=l_m=0$, and we have
\begin{eqnarray}
\label{litthe}
\sum_{E\in E_n^{0,\ldots,0}}
B_{a^r,a^m} \,
B_{\nu-E}
 =  B_{a^{r+m}}\, B_{\nu}=
\!\! \!\!
\sum_{L=(l_1,\ldots,l_r)\in[m]}
\sum_{E\in E_m^{L}}
(-t)^{|L|}
B_{a+l_1,\cdots,a+l_r} \,
B_{a^m-E,\nu}
\,.
\end{eqnarray}
Moreover, since $E^{L}_m =E^{\sigma L}_m$ for any permutation
$\sigma$, (\ref{litthe}) can be written as
\begin{equation}
B_{a^{r+m}}\,B_{\nu}
= \sum_{\lambda\in\P^r_{\leq m}} (-t)^{|\lambda|}
\left( \sum_{\sigma \in S_r;\sigma(\lambda)~{\rm{ distinct}}}
B_{\sigma(\lambda)+a^r} \right)
\left( \sum_{E \in E^{\lambda}_m} B_{a^m-E,\nu}
\right) \, .
\end{equation}
We can now use Lemmas~\ref{iden1} and \ref{iden2} to obtain 
\begin{equation} \label{lasteq}
\begin{split} 
B_{a^{r+m}} \,  B_{\nu} & = \sum_{\lambda \in \P^r_{\leq m}} (-t)^{|\lambda|} 
 \left(\sum_{\mu \in \P^r} B_{\mu+a^r} K_{\lambda \mu}^{-1} \right)
\left(\sum_{\omega} \sum_{\gamma \in \P^m} \sum_{\rho \in \P^m_{\leq a}}
 K_{\gamma \rho}^{-1}
K_{\omega \lambda} K_{\omega' \gamma} \, B_{a^m-\rho^R,\nu} \right)\\
& = 
\sum_{\rho \in \P^m_{\leq a}}
\sum_{\gamma \in \P^m} \sum_{\mu \in \P^r} (-t)^{|\mu|}
 B_{\mu+a^r}  
\sum_{\omega} 
\left( \sum_{\lambda \in \P^r_{\leq m} }
K_{\omega \lambda} K_{\lambda \mu}^{-1} \right)   K_{\gamma \rho}^{-1}
 K_{\omega' \gamma} \, B_{a^m-\rho^R,\nu} \\
\end{split}
\end{equation}
It happens that we can sum over all $\lambda$ since any term
$\lambda \not \in \P^{r}_{\leq m}$ vanishes by the known property
$K_{\lambda \mu}=K_{\lambda \mu}^{-1}=0$ 
if $\lambda \not \geq \mu$.  
That is, if $\lambda \not \in \P^r$, 
$K_{\lambda \mu}^{-1}=0$ since $\mu \in \P^r$. 
Moreover, if $\lambda_1>m$, $K_{\omega \lambda}=0$ 
if $\omega_1 \not > m$.  But if $\omega_1 > m$, we have 
$\omega' \not \in \P^m$, and thus $K_{\omega' \gamma}=0$ 
since $\gamma \in \P^m$.  Therefore,
\begin{equation} 
\begin{split} 
B_{a^{r+m}} \,  B_{\nu}  
& =
\sum_{\rho \in \P^m_{\leq a}} 
\sum_{\gamma \in \P^m} \sum_{\mu \in \P^r} (-t)^{|\mu|}
B_{\mu+a^r} 
\sum_{\omega} 
 \left( \sum_{\lambda  }
K_{\omega \lambda} K_{\lambda \mu}^{-1} \right)   K_{\gamma \rho}^{-1}
 K_{\omega' \gamma} \, B_{a^m-\rho^R,\nu} \\
& = \sum_{\rho \in \P^m_{\leq a}} \sum_{\mu \in \P^r} (-t)^{|\mu|}
 B_{\mu+a^r} \left( \sum_{\gamma \in \P^m}   
K_{\mu' \gamma}  K_{\gamma \rho}^{-1}
 \right) \, B_{a^m-\rho^R,\nu}  \, ,
\end{split}
\end{equation}
since $\sum_{\omega} \delta_{\omega \mu} K_{\omega' \gamma}=K_{\mu' \gamma}$.
Again, if $\gamma \not \in \P^m$, we have $K_{\gamma \rho}^{-1}=0$ since
$\rho \in \P^m_{\leq a}$.  Thus, we have
\begin{equation}
\begin{split}
B_{a^{r+m}} \,  B_{\nu}  
& = \sum_{\rho \in \P^m_{\leq a}}
 \sum_{\mu \in \P^r} (-t)^{|\mu|}
B_{\mu+a^r} \left( \sum_{\gamma }   K_{\mu' \gamma}  K_{\gamma \rho}^{-1}
\right) \, B_{a^m-\rho^R,\nu}  \\
& = \sum_{\rho \in \P^m_{\leq a}}
 \sum_{ \mu \in \P^r}  (-t)^{|\rho|} \, \delta_{\mu' \rho}
\,
 B_{\mu+a^r}  \, B_{a^m-\rho^R,\nu} \\
& = \sum_{ \mu \in \P^{r}_{\leq m}}  (-t)^{|\mu|} \,
 B_{\mu+a^r}  \, B_{a^m-(\mu')^R,\nu} \, ,
\end{split}
\end{equation}
with the restriction that the summand vanishes if $(a^m-(\mu')^R)$
is not a partition.  \hfill $\square$

\section{$k$-split polynomials} \label{sec3}
Recent developments in the study of the symmetric function
space have centered around a filtration,  
$\Lambda^{(1)}_t \subseteq \Lambda^{(2)}_t \subseteq \cdots 
\subseteq \Lambda^{(\infty)}_t=\Lambda$, 
given by the subspaces
\begin{equation} \label{eqlinearspan}
\Lambda^{(k)}_t = 
{\mathcal L} \left\{ H_\lambda[X;t] \right\}_{\lambda;\lambda_1\leq k} 
=
{\mathcal L} \left\{ H_\lambda[X;q,t] \right\}_{\lambda;\lambda_1\leq k} \, .
\end{equation}
This filtration provides a convenient environment for the
generalization of the theory of Schur functions,
and for the study of Macdonald polynomials \cite{[LLM],[LM2]}.
Several new families of polynomials were introduced in \cite{[LM2]}
as an approach to studying the spaces $\Lambda^{(k)}_t$.
In this section, we use results from Section~\ref{SecVer} 
to prove properties related to one of these families; 
the $k$-split polynomials.

\begin{definition}
\label{Defksp}
For a $k$-bounded partition $\lambda$,
let $\lambda^{\to k}=(\lambda^{(1)},\lambda^{(2)},\dots)$.
The $k$-split polynomials are defined recursively by
\begin{equation} \label{recursiB}
G_{\lambda}^{(k)}[X;t] = 
B_{\lambda^{(1)}} {G}^{(k)}_{(\lambda^{(2)},\lambda^{(3)},\dots)} [X;t] \, ,
\quad\text{with}\quad G^{(k)}_{()}=1
\, .
\end{equation}
\label{defksp}
\end{definition}

It was shown in \cite{[LM2]} that
the $k$-split polynomials form a basis for $\Lambda_t^{(k)}$, and that
vertex operators indexed by partitions with 
hook-length not larger than $k$ leave this space 
invariant.  More precisely, 

\begin{property}
\cite{[LM2]}
\label{lempreserve}
If $\lambda$ is a partition with $h_M(\lambda)\leq k$,
then $B_\lambda\, f\in\Lambda_t^{(k)}$ for any
$f\in\Lambda_t^{(k)}$. 
\end{property}

\noindent
Further, it was shown that $B_i$ acts invariantly
on subspaces of $\Lambda^{(k)}_t$, defined 
for integers $a \leq  k$ by
\begin{equation}
\Lambda^{(a,k)}_t = {\mathcal L}\{ H_\lambda[X;t] \}_{a\leq \lambda_1\leq k}
\, . 
\end{equation}
Note that in the case $a\leq 0$, we simply have 
$\Lambda^{(a,k)}_t=\Lambda^{(k)}_t$, and that $\Lambda^{(a,k)}_t$
can also be defined as
\begin{equation}
\Lambda^{(a,k)}_t = {\mathcal L}\{ G_\lambda[X;t] \}_{a\leq \lambda_1\leq k}
\, ,
\end{equation}
since the transition matrix between the
two bases is upper triangular \cite{[LM2]}.  

\begin{property}
\cite{[LM2]}
\label{j1}
If $i$ is an integer such that $i \leq  k$, then 
$B_i\,f\in\Lambda_t^{(i,k)} \subseteq \Lambda^{(k)}_t$ for all
$f\in \Lambda_t^{(k)}$.
\end{property}

We start by extending these properties and then discover,
more specifically, that there exist subspaces of $\Lambda^{(a,k)}_t$
that are invariant under a special set of
the $B_\lambda$ operators.
These subspaces are defined, for nonnegative integers $j \leq k$, 
\begin{equation}
\Omega^{(k)}_j = {\mathcal L} 
\bigl \{ G_{\lambda}^{(k)}[X;t] \bigr \}_{\lambda_1=j}\, .
\end{equation}
Our generalization of Property~\ref{lempreserve} is

\bigskip

\begin{property} 
\label{lemigen}
If $\lambda$ is a partition with $h_M(\lambda)\leq k$, then
$B_\lambda\,f\in\Lambda_t^{(\lambda_1,k)}$ for all
$f\in \Lambda_t^{(k)}$.
\end{property}
\noindent{\bf Proof.}\quad
Definition \eqref{jtt} gives
\begin{equation}
\begin{split}
B_\lambda
& =
\prod_{2\leq j \leq \ell(\lambda)} (1-te_{1j}) B_{\lambda_1}
\prod_{2\leq i <j \leq \ell(\lambda)} (1-te_{ij}) B_{\lambda_2} \cdots
B_{\lambda_{\ell(\lambda)}}
\\
&=
\prod_{2\leq j \leq \ell(\lambda)} (1-te_{1j}) B_{\lambda_1}
B_{\hat\lambda}
\end{split}
\end{equation}
Since $e_{1j}$ increases the index of $B_{\lambda_1}$,
we have
\begin{equation}
B_{\lambda}= 
\sum_{i=0}^{\ell(\lambda)-1}
c_i(t) \, B_{\lambda_1+i} \, O_i \, ,\quad \qquad c_i(t) \in \mathbb Z[t] \, ,
\label{te}
\end{equation}
where $O_i$ is a product of $B_{j}$'s with $j\leq k$.
Let $f\in\Lambda_t^{(k)}$ and note that Property~\ref{j1} 
implies $B_j\,f\in\Lambda_t^{(k)}$ for $j\leq k$. 
Therefore, $O_i \cdot f \in \Lambda^{(k)}_t$
and again by Property~\ref{j1},
$B_{\lambda_1+i}\, O_i \cdot f \in \lambda^{(\lambda_1+i,k)}_t
\subseteq\Lambda_t^{(\lambda_1,k)}$
since $\lambda_1+i\leq \lambda_1+\ell(\lambda)-1=h_M(\lambda)\leq k$.
\hfill $\endprf$

\bigskip

\begin{property}
\label{lemikplus}
If $i$ is an integer such that $1\leq i<k$ then 
$B_i\,f\in\Lambda_t^{(i+1,k)}$ for all
$f\in \Lambda_t^{(i+1,k)}$.
\end{property}
\noindent{\bf Proof.}\quad
\quad 
Let $f \in \Lambda^{(i+1,k)}_t$ and assume without
loss of generality that $f= H_{\lambda}[X;t]=B_{\lambda_1}
H_{\hat\lambda}[X;t]$ for $\lambda$ with $\lambda_1 > i$.   
When $\lambda_1=i+1$, we have 
$B_i f = B_{i}B_{i+1}H_{\hat \lambda}[X;t]
=B_{i+1} B_{i} H_{\hat \lambda}[X;t]$
by the commutation relation \eqref{commuB}.
Property~\ref{j1} then implies that
$B_if=B_{i+1}B_{i} H_{\hat \lambda}[X;t] \in \Lambda^{(i+1,k)}_t$.
In the case that $k \geq \lambda_1> i+1$, again by relation
\eqref{commuB}, we have
\begin{equation}
B_i \, H_{\lambda}[X;t]= B_{i} \, B_{\lambda_1} \, H_{\hat \lambda}[X;t] = 
\bigl( t B_{\lambda_1} B_i  + 
t B_{i+1} B_{\lambda_1-1}
- B_{\lambda_1-1} B_{i+1} \bigr) \cdot H_{\hat \lambda} \, .
\end{equation} 
The three terms in the right hand side are all of the type
$B_{a} B_{b}H_{\hat \lambda}$, where $a > i $ and $b \leq k$.
Therefore, 
$B_{a}B_{b}H_{\hat \lambda} \in \Lambda^{(i+1,k)}_t$
again by Property~\ref{j1}.
\hfill $\square$

We now prove that the subspaces $\Omega_j^{(k)}$
are invariant under the action of operators indexed by $k$-rectangles,
using the following lemma:

\begin{lemma} \label{lemomega}
If $\lambda$ is a partition with $h_M(\lambda)=k$ and
$\lambda_L \geq j$
then $B_{\lambda} \, f \in \Omega^{(k)}_{\lambda_1}$
for all $f \in \Omega^{(k)}_j$.
\end{lemma}
\noindent {\bf Proof.} \quad  
Letting $f=G_{\mu}^{(k)}[X;t]$ with $\mu_1=j$,
we have that $(\lambda,\mu)$ is a partition since
$\mu_1\leq\lambda_L$.  
If
$\mu^{\to k}= (\mu^{(1)}, \mu^{(2)},\dots)$,
we thus have
$\big(\lambda,\mu\big)^{\to k}=(\lambda,\mu^{(1)},\mu^{(2)},\dots)$
since $h_M(\lambda)=k$.  Therefore, 
\begin{equation}
B_{\lambda} \, f = B_{\lambda} \left( 
B_{\mu^{(1)}} \, B_{\mu^{(2)}}\cdots\right) \cdot 1 =
G^{(k)}_{(\lambda,\mu)}[X;t] 
\end{equation}
by definition,
and thus $B_{\lambda} \, f \in \Omega^{(k)}_{\lambda_1}$.
\hfill $\square$

We are now prepared to prove the final result of this section.

\bigskip

\begin{theorem} \label{propnouv}
If $j\leq k$ is a nonnegative integer
then
$B_{\ell^{k+1-\ell}} \, f \in \Omega^{(k)}_{\max(j,\ell)} $ 
for all $f \in \Omega^{(k)}_j$. 
\end{theorem}
\noindent {\bf Proof.} \quad 
We have either (a) $\max(j,\ell)=\ell$ or (b) $\max(j,\ell)=j$.
In case (a), letting $\lambda=(\ell^{k+1-\ell})$ in 
Lemma~\ref{lemomega}, the assertion holds.  
Thus assume $j>\ell$ and
let $f= G_{\nu}^{(k)}[X;t]$ where $\nu_1=j$.
If $\nu^{\to k}=(\nu^{(1)},\nu^{(2)},\ldots)$
then there are three possibilities:
$(i)$ $h_M(\nu^{(1)})<k$,
$(ii)$ $h_M(\nu^{(1)})=k$
and $\nu^{(1)}_{L} < \ell$, or $(iii)$  $h_M(\nu^{(1)})=k$
and $\nu^{(1)}_{L} \geq \ell$, where $L=\ell(\nu^{(1)})$. 

In case $(i)$, $h_M(\nu^{(1)})< k$ implies $\nu^{\to k}= 
(\nu^{(1)})$ and thus $f=B_{\nu^{(1)}} \cdot 1$.
Identity~\ref{identity4} then gives
\begin{equation}
B_{\ell^{k+1-\ell}} \, B_{\nu^{(1)}} \cdot 1 = 
\sum_{i} (-t)^{c(i)} B_{\rho^{(i)}} \, B_{\gamma^{(i)}} \cdot 1 \, ,
\end{equation}
where the partitions $\rho^{(i)}$ and $\gamma^{(i)}$ 
are such that $(\rho^{(i)},\gamma^{(i)})^{\to k}=
(\rho^{(i)},\gamma^{(i)})$  and $\rho^{(i)}_1=j$. Therefore,
\begin{equation}
B_{\ell^{k+1-\ell}} \, B_{\nu^{(1)}} \cdot 1 = 
\sum_{i} (-t)^{c(i)} G^{(k)}_{(\rho^{(i)}, \gamma^{(i)})}[X;t] 
\in \Omega^{(k)}_j \, .
\end{equation}

For case $(ii)$, 
we have $f= B_{\nu^{(1)}} B_{\nu^{(2)}} \cdots 1$.
Identity~\ref{identity2} states that
\begin{equation}
B_{\ell^{k+1-\ell}} \, B_{\nu^{(1)}}  = 
\sum_{i} (-t)^{c(i)} B_{\rho^{(i)}} \, B_{\gamma^{(i)}} \,  ,
\end{equation}
where $\rho^{(i)}$ and $\gamma^{(i)}$ are such that
$(\rho^{(i)},\gamma^{(i)},\nu^{(2)},\dots)^{\to k}=
(\rho^{(i)},\gamma^{(i)},\nu^{(2)},\dots)$ 
and $\rho^{(i)}_1=j$.  Thus
\begin{equation}
B_{\ell^{k+1-\ell}} \, f = 
\sum_{i} (-t)^{c(i)} G_{(\rho^{(i)}, \gamma^{(i)},
\nu^{(2)},\ldots)}[X;t] \in \Omega^{(k)}_j \, .
\end{equation}

Finally in case $(iii)$, 
first $t$-commute $B_{\ell^{k+1-\ell}}$ 
with the operators $B_{\nu^{(i)}}$,
using Identity~\ref{identity3}, until
\begin{equation} \label{eqprod}
B_{\ell^{k+1-\ell}} \,  B_{\nu^{(1)}} 
\cdots  B_{\nu^{(m)}} \, B_{\nu^{(m+1)}} \cdots 1 = 
t^{*} B_{\nu^{(1)}} \cdots B_{\nu^{(m)}}\, B_{\ell^{k+1-\ell}} \,
 B_{\nu^{(m+1)}} \cdots 1 \, ,
\end{equation}
where $*$ is a power of $t$, and where
$\nu^{(m)}$ is such that $\nu^{(m)}_L \geq \ell$,
while $\nu^{(m+1)}_1 < \ell$, $h_M(\nu^{(m+1)})<k$ 
or $\nu^{(m+1)}_L <\ell$.  
In any of these scenarios, if $\mu^{\to k}=(\nu^{(m+1)},\dots)$
then $\mu$ satisfies the conditions of cases 
$(a)$, $(i)$ or $(ii)$, resp.
We thus have
\begin{equation}
f=
B_{\ell^{k+1-\ell}}  B_{\nu^{(m+1)}} \cdots 1 = B_{\ell^{k+1-\ell}}\, 
G_{\mu}^{(k)}[X;t] 
\; \in\; \Omega^{(k)}_{\max \{\ell,\mu_1=\nu^{(m+1)}_1 \}}\, .
\end{equation}
Since $\max\{\ell,\nu^{(m+1)}_1\}\leq\nu^{(m)}_L$, 
applying $B_{\nu^{(m)}}$ to $f$
then produces an element of $\Omega^{(k)}_{\nu^{(m)}_1}$
by Lemma~\ref{lemomega}. 
Applying $B_{\nu^{(1)}},\ldots,B_{\nu^{(m-1)}}$ by the
same argument gives
$B_{\nu^{(1)}}\cdots B_{\nu^{(m)}}\, B_{\ell^{k+1-\ell}} \,
 B_{\nu^{(m+1)}} \cdots 1\in \Omega^{(k)}_{\nu^{(1)}_1}=
\Omega^{(k)}_{\nu}$ proves the theorem from (\ref{eqprod}).
\hfill $\square$

\section{$k$-Schur functions} \label{sec4}

The $k$-split polynomials 
play a crucial role in the generalization of
the theory of symmetric functions since they 
are essential for the construction of another family 
of polynomials called the $k$-Schur functions. 
The characterization of this family relies on 
a projection operator that acts linearly on 
$\Lambda^{(k)}_t$, for nonnegative integers 
$j \leq k$, by
\begin{equation}
T_j^{(k)} \, G_\lambda^{(k)}[X;t] = 
\begin{cases}
G_\lambda^{(k)}[X;t] & \text{if } \lambda_1 = j\\
0 & \text{otherwise }
\end{cases}
\, .
\end{equation}

\begin{definition}  
\label{defks}
Let $\lambda$ be a $k$-bounded partition.  
The $k$-Schur functions are defined recursively
by, 
\begin{equation}
s_{\lambda}^{(k)}[X;t] = {T}_{\lambda_1}^{(k)} B_{\lambda_1} 
s_{(\lambda_1,\lambda_2,\dots)}^{(k)}[X;t] \, ,
\quad\text{where}\quad
s_{()}^{(k)}[X;t]=1\, .
\end{equation}
\end{definition}

The $k$-Schur functions are believed to play 
a role in $\Lambda^{(k)}_t$ that is analogous
to the role the Schur functions have in $\Lambda$.
It was shown that these functions
form  a basis for $\Lambda^{(k)}_t$ and that they 
reduce to the Schur functions themselves when $k\to\infty$. 
That is, $s_{\lambda}^{(k)}[X;t] = s_{\lambda}[X]$
when $h_M(\lambda)\leq k$.
Several other properties supporting the claim that the $k$-Schur functions 
generalize the theory of Schur functions are given in \cite{[LLM],[LM2]}.

A continuation of the study in Section~\ref{sec3} 
regarding the action of operators indexed by $k$-rectangles
led to an observation that not all of the 
$k$-Schur functions
need to be constructed using Definition~\ref{defks}.  
We found that for each $k$, there is a subset of $s_\lambda^{(k)}[X;t]$
called the irreducible $k$-Schur functions, from which all other 
$k$-Schur functions may be constructed by simply applying a 
succession of operators indexed by $k$-rectangles.
This subset consists of the special set of $k$-Schur functions
indexed by irreducible partitions;
$k$-bounded partitions with no 
more than $i$ parts equal to $k-i$, for $i=0,\dots,k-1$.  

\begin{definition}
A $k$-Schur function indexed by an irreducible partition 
is said to be a irreducible.
Otherwise, the $k$-Schur function is called reducible.
\label{defindecom}
\end{definition}

\noindent For example, the irreducible $k$-Schur functions for $k=1,2,3$ are
\begin{equation}
\begin{split}
k=1 \, &:  \qquad s_0^{(1)} \, ,\\
k=2 \, &: \qquad s_0^{(2)} \, , \quad s_1^{(2)} \, ,\\
k=3 \, &: \qquad s_0^{(3)} \, , \quad s_1^{(3)}\, ,
\quad s_{2}^{(3)}\, , \quad s_{1,1}^{(3)} \, , \quad s_{2,1}^{(3)} 
\, , \quad s_{2,1,1}^{(3)}\, .
\end{split} \end{equation}

\noindent These examples support the following property;
\begin{property}  \cite{[LLM]}
There are $k!$ distinct $k$-irreducible partitions.
\label{lemmakfac}
\end{property}

The main result of this section is to prove that an operator 
indexed by a $k$-rectangle $R$ sends, up to a constant, $s_\lambda^{(k)}[X;t]$ 
to $s_{R\cup\lambda}^{(k)}[X;t]$.  Consequently, given the 
$k!$ irreducible $k$-Schur functions,  any reducible $k$-Schur 
function can be obtained by applying a sequence of 
operators indexed by $k$-rectangles to the appropriate 
irreducible $k$-Schur function.  
We first give several properties of the
operators $B_{\lambda}$ and $T_j^{(k)}$
and then, using our results from Section~\ref{sec3},
we will prove our main result.

\begin{property} \label{lemsepara}
\cite{[LM2]}
If $\lambda$ is a partition with $h_M(\lambda) \leq k$
then
${T}^{(k)}_{\lambda_1} B_{\lambda}  f = 
{T}^{(k)}_{\lambda_1} B_{\lambda_1} B_{\hat \lambda} f $
for all $f \in \Lambda^{(k)}_t$.
\end{property}

\bigskip

\begin{property}  \label{lemcommu}
If $k\!\geq j\! >\! \ell$ are nonnegative integers,
$T_{j}^{(k)} B_{\ell^{k+1-\ell}}\, f= 
B_{\ell^{k+1-\ell}} T_{j}^{(k)} \, f$
for all $f \in \Lambda^{(k)}_t$.
\end{property}
\noindent {\bf Proof.} \quad 
It suffices to consider $f= G_{\lambda}^{(k)}[X;t]$.  
By the definition of ${T}^{(k)}_j$, we have
\begin{equation} \label{eqcom1}
B_{\ell^{k+1-\ell}} \, {T}_{j}^{(k)} \, f =
\begin{cases}
0 & {\rm if~} j \neq \lambda_1 \\
B_{\ell^{k+1-\ell}} \, f  & {\rm if~} j = \lambda_1 \\
\end{cases} \, .
\end{equation}
On the other hand, consider
$T_j^{(k)} B_{\ell^{k+1-\ell}}f$.
By Theorem~\ref{propnouv},
$B_{\ell^{k+1-\ell}}\,f\in \Omega_{\max\{\ell,\lambda_1\}}^{(k)}$.
If $j=\lambda_1>\ell$, then
$B_{\ell^{k+1-\ell}}\,f\in \Omega_{j}^{(k)}\implies
T_j^{(k)} B_{\ell^{k+1-\ell}}\,f = B_{\ell^{k+1-\ell}}\,f$.
If $j\neq\lambda_1$ then either
$B_{\ell^{k+1-\ell}}\,f\in \Omega_{\ell<j}^{(k)}$ or
$B_{\ell^{k+1-\ell}}\,f\in \Omega_{\lambda_1\neq j}^{(k)}$.
Both cases vanish under the action of $T_j^{(k)}$ and thus
we prove our claim. \hfill $\square$

We can now show that an operator indexed by a $k$-rectangle 
acts simply on a $k$-Schur function.

\bigskip

\begin{theorem} \label{theorec}
If $\mu,\nu,\lambda$
are partitions where $\lambda=(\mu,\nu)$ and
$\mu_L>\ell\geq \nu_1$, then
\begin{equation}
B_{\ell^{k+1-\ell}} \, s_{\lambda}^{(k)}[X;t] = t^{|\mu|-\ell(\mu) \ell}
s_{(\ell^{k+1-\ell})\cup\lambda}^{(k)}[X;t] \, .
\end{equation}
\end{theorem}
\noindent {\bf Proof.} \quad  Let $\ell(\mu)=M$. 
Since  $\lambda_M=\mu_L>\ell$ and
$\lambda_{M+1}=\nu_1\leq\ell$, we have 
\begin{equation}
\label{shh}
B_{\ell^{k+1-\ell}} \, s_{\lambda}^{(k)}[X;t] = 
t^{\lambda_1+\dots+\lambda_M-M\ell}
T^{(k)}_{\lambda_1}
B_{\lambda_1} 
T^{(k)}_{\lambda_2}
B_{\lambda_2} 
\dots
T^{(k)}_{\lambda_{M}}
B_{\lambda_{M}} 
B_{\ell^{k+1-\ell}}
\, s_{\nu}^{(k)}[X;t]
\, ,
\end{equation}
by iterating the following argument $M$ times:
If $\hat\lambda=(\lambda_2,\ldots,\lambda_M)$,
by definition of $s_\lambda^{(k)}$ we have
\begin{equation}
B_{\ell^{k+1-\ell}} \, s_{\lambda}^{(k)}[X;t] = 
B_{\ell^{k+1-\ell}} \, 
T^{(k)}_{\lambda_1} \, B_{\lambda_1} \,  s_{\hat \lambda}^{(k)}[X;t] \, .
\end{equation}
Since $s_{\hat \lambda}^{(k)}\in\Lambda^{(k)}_t$,
Property~\ref{lempreserve}
gives that $B_{\lambda_1} s_{\hat \lambda}^{(k)}\in\Lambda^{(k)}_t$.
Property~\ref{lemcommu} then implies that
$T_{\lambda_1}^{(k)}$ commutes with $B_{\ell^{k-\ell+1}}$ since
$\lambda_1>\ell$.  We thus have
\begin{equation}
B_{\ell^{k+1-\ell}} \, s_{\lambda}^{(k)}[X;t] = 
T^{(k)}_{\lambda_1} B_{\ell^{k+1-\ell}}
B_{\lambda_1} \,  s_{\hat \lambda}^{(k)}[X;t] 
\, .
\end{equation}
Furthermore, $B_{\ell^{k-\ell+1}}$ $t$-commutes with 
$B_{\lambda_1}$ by Identity~\ref{identityzero}.  Therefore,
\begin{equation}
B_{\ell^{k+1-\ell}} \, s_{\lambda}^{(k)}[X;t] = 
t^{\lambda_1-\ell}
T^{(k)}_{\lambda_1}
B_{\lambda_1} 
B_{\ell^{k+1-\ell}}
\,  s_{\hat \lambda}^{(k)}[X;t] 
\, .
\end{equation}

Now if we can show, for a partition $\nu$ with $\nu_1\leq \ell$, 
that
\begin{equation}
\label{lastlem}
B_{\ell^{k+1-\ell}} \,  s_{\nu}^{(k)}[X;t] 
= 
\left( T_\ell^{(k)}\, B_{\ell}\right)^{k+1-\ell} \,  s_{\nu}^{(k)}[X;t] 
\, ,
\end{equation}
then putting this into \eqref{shh} 
implies our result by definition, 
since $\lambda_M > \ell\geq\lambda_{M+1}$.
To prove \eqref{lastlem}, we have
$s_{\nu}^{(k)}\in \Omega^{(k)}_{\nu_1}$ by definition,
which implies
$B_{\ell^{k+1-\ell}} s_{\nu}^{(k)}\in\Omega^{(k)}_{\ell}$
by Theorem~\ref{propnouv} since $\nu_1\leq \ell$. 
Therefore, $B_{\ell^{k+1-\ell}} s_{\nu}^{(k)}$
is invariant under $T_\ell^{(k)}$
and \eqref{lastlem} is equivalent to
\begin{equation}
{T}^{(k)}_{\ell} B_{\ell^{k+1-\ell}}s_{\nu}^{(k)}[X;t]
= 
\left( T_\ell^{(k)}\, B_{\ell}\right)^{k+1-\ell} \,  s_{\nu}^{(k)}[X;t] 
\, .
\end{equation}
Using Property~\ref{lemsepara}, the right hand side may be
written
\begin{equation}
\label{pe}
{ T}^{(k)}_{\ell}\, 
B_{\ell^{k+1-\ell}} \, s_{\nu}^{(k)}[X;t] = { T}^{(k)}_{\ell}\, 
B_{\ell} \, B_{\ell^{k-\ell}} \, s_{\nu}^{(k)}[X;t] \, . 
\end{equation}
Now, 
$B_{\ell^{k-\ell}} s_{\nu}^{(k)}[X;t]\in\Lambda^{(\ell,k)}_t$
by Property~\ref{lemigen}.
Further, since any element of $\Lambda^{(\ell,k)}_t$ 
can be decomposed into the sum of two functions,
$f\in\Omega^{(k)}_\ell$ and $g\in\Lambda^{(\ell+1,k)}_t$,
we have
$B_{\ell^{k-\ell}} \, s_{\nu}^{(k)}[X;t] = f+g$.
Therefore,
\begin{equation}
{ T}^{(k)}_{\ell}\, B_{\ell} \, (f+g)=
{ T}^{(k)}_{\ell}\, B_{\ell}\,  { T}^{(k)}_{\ell}\, f
=
{ T}^{(k)}_{\ell}\, B_{\ell}\,  { T}^{(k)}_{\ell}
\,( f+g)\,
\end{equation}
since $B_{\ell}\cdot g \in \Lambda^{(\ell+1,k)}_t$
by Property~\ref{lemikplus} implies 
$T_\ell^{(k)}B_{\ell} \cdot g=0$.
Thus, \eqref{pe} gives
\begin{equation} 
T^{(k)}_\ell\, B_{\ell^{k+1-\ell}} \, s_{\nu}^{(k)}[X;t] =
{ T}^{(k)}_{\ell}B_{\ell} \, { T}^{(k)}_{\ell}
B_{\ell^{k-\ell}} \, s_{\nu}^{(k)}[X;t] \, .
\end{equation}
Repeating this argument $k-\ell$ times proves (\ref{lastlem}).
\hfill $\square$

\medskip

When $t=1$, our theorem reduces to the simple expression:

\medskip

\begin{corollary} 
\label{conjrecschur} 
If $\lambda$ is a $k$-bounded partition, then
\begin{equation}
s_{\ell^{k+1-\ell}}[X] \, s_{\lambda}^{(k)}[X] = 
s_{(\ell^{k+1-\ell})\cup\lambda}^{(k)}[X] \, .
\end{equation}
\end{corollary}

\noindent
The role of Schur functions indexed by $k$-rectangles
in the subring $\Lambda^{(k)}$ leads naturally
to the study of  the quotient ring $\Lambda^{(k)}/{\mathcal I}_k$,
where ${\mathcal I}_k$ denotes the ideal generated by
$s_{(\ell^{k+1-\ell})}[X]$.
It is known that the dimension of this
quotient ring is $k!$ since

\begin{proposition}
\cite{[LLM]}
The homogeneous functions indexed by $k$-irreducible partitions
form a basis of the quotient ring $\Lambda^{(k)}/{\mathcal I}_k$.
\end{proposition}

\noindent
Corollary~\ref{conjrecschur} then implies

\medskip

\begin{theorem}
The irreducible $k$-Schur functions form a basis of the quotient ring
$\Lambda^{(k)}/{\mathcal I}_k$.
\end{theorem}
The irreducible $k$-Schur function basis thus offers a simple
method of performing operations in this quotient ring:  
first work in $\Lambda^{(k)}$ using $k$-Schur functions and then
replace by zero all the $k$-Schur functions indexed by partitions which are
not $k$-irreducible.

\section{Appendix}
\label{A1}
\begin{formula} \label{kostkaplus}
 Let $\mu,\lambda \in \P^m$.   
The Kostka matrix is such that
\begin{equation}
K_{\lambda \mu}^{-1} = K_{\lambda+a^m,\mu+a^m}^{-1} \, . 
\end{equation}
\end{formula}  
\noindent{\bf{Proof.}}\quad  
We have
\begin{equation}\label{eqrec1} 
\begin{split}
m_{\lambda+a^m}[x_1+\cdots+x_m] & = (x_1\cdots x_m)^a \, 
m_{\lambda}\left[{x_1} + \cdots + {x_m} \right] \\
& = (x_1 \cdots x_m)^a \, \sum_{\mu} K_{\lambda \mu}^{-1} \,
s_{\mu}\left[ {x_1} + \cdots +{x_m} \right] \\
& =  \sum_{\mu} K_{\lambda \mu}^{-1} \,
s_{\mu+a^m}\left[ {x_1} + \cdots + {x_m} \right] \, ,
\end{split}
\end{equation}
since $(x_1 \cdots x_m)^a \, 
s_{\mu}\left[ {x_1} + \cdots +{x_m} \right]=
s_{\mu+a^m}\left[ {x_1} + \cdots +{x_m} \right]$.
By definition
\begin{equation} \label{eqnu2}
m_{\lambda+a^m}[x_1+\cdots+x_m] = 
\sum_{\nu} K_{\lambda+a^m,\nu}^{-1} \,
s_{\nu}\left[ {x_1} + \cdots + {x_m} \right] \, . 
\end{equation}
Taking the coefficient of 
$s_{\mu+a^m}$
in (\ref{eqrec1}) and (\ref{eqnu2}), we get 
$K_{\lambda \mu}^{-1} = K_{\lambda+a^m,\mu+a^m}^{-1}$, as claimed.
\hfill $\square$

\begin{formula} \label{formu1}
 Let $\lambda$ be a partition of at most $m$ parts, 
and let $a \geq \lambda_1$.  Then
\begin{equation}
s_{\lambda}\left[ \frac{1}{x_1} + \cdots + \frac{1}{x_m} \right](x_1 \cdots x_m)^a
= s_{a^m-\lambda^R}[x_1+\cdots+x_m] \, .
\end{equation}
\end{formula}
\noindent{\bf{Proof.}}\quad  In \cite{[Lint]} ($Sf5$), one finds
\begin{equation}
\begin{split}
s_{\lambda}\left[\frac{1}{x_1}+ \cdots+\frac{1}{x_m}\right] & = 
s_{(a^m/\lambda')'}[x_1+\cdots+x_m]/s_{1^m}[x_1+\cdots+x_m]^a \\
& = s_{(a^m/\lambda')'}[x_1+\cdots+x_m]/(x_1\cdots x_m)^a  \, .
\end{split}
\end{equation}
Following from the Jacobi-Trudi determinantal expression for skew Schur
functions, we have
\begin{equation}  
\begin{split}
s_{\mu/\lambda} = \det \left|s_{\mu_i-\lambda_j-i+j} \right|_{1\leq i,j
\leq m}
& = \det \left|s_{\mu_{n+1-j}-\lambda_{n+1-i}-(n+1-j)+(n+1-i)} 
\right|_{1\leq i,j \leq m} \\
& =  \det \left|s_{\mu_{n+1-j}-\lambda_{n+1-i}-i+j} 
\right|_{1\leq i,j \leq m} \\
\end{split} \, .
\end{equation}
Therefore, if $\mu=(a^m)$, we obtain
\begin{equation}
s_{a^m/\lambda} = \det \left|s_{a-\lambda_{n+1-i}-i+j} 
\right|_{1\leq i,j \leq m} = s_{a^m-\lambda^R} \, ,
\end{equation}
completing the proof.  \hfill$\square$

\begin{formula} \label{kostkamoins}
 Let $\mu,\lambda \in \P^m$,
and let $a$ be an integer such that $a \geq \lambda_1$ and
$a \geq \mu_1$.  Then, the Kostka matrix is such that
\begin{equation}
K_{\lambda \mu}^{-1} = K_{a^m-\lambda^R,a^m-\mu^R}^{-1} \, . 
\end{equation}
\end{formula}  
\noindent{\bf{Proof.}}\quad  
Using Formula~\ref{formu1}, we have
\begin{equation} \label{eqrec}
\begin{split}
m_{a^m-\lambda^R}[x_1+\cdots+x_m] & = (x_1\cdots x_m)^a \, 
m_{\lambda}\left[ \frac{1}{x_1} + \cdots + \frac{1}{x_m} \right] \\
& = (x_1 \cdots x_m)^a \, \sum_{\mu} K_{\lambda \mu}^{-1} \,
s_{\mu}\left[ \frac{1}{x_1} + \cdots + \frac{1}{x_m} \right] \\
& =  \sum_{\mu} K_{\lambda \mu}^{-1} \,
s_{a^m-\mu^R}\left[ {x_1} + \cdots + {x_m} \right] 
\end{split}
\end{equation}
Since, by definition,
\begin{equation} \label{eqnu}
m_{a^m-\lambda^R}[x_1+\cdots+x_m] = 
\sum_{\nu} K_{a^m-\lambda^R,\nu}^{-1} \,
s_{\nu}\left[ {x_1} + \cdots + {x_m} \right] \, , 
\end{equation}
taking the coefficient of $s_{a^m-\mu^R}$
in (\ref{eqrec}) and (\ref{eqnu}), we get 
$K_{\lambda \mu}^{-1} = K_{a^m-\lambda^R,a^m-\mu^R}^{-1}$, as claimed.
\hfill $\square$

\begin{formula} \label{monoschur}
If $\lambda=(\lambda_1,\dots,\lambda_n)$, then  
\begin{equation}
m_{\lambda}= \sum_{\sigma \in S_n;
\sigma(\lambda) {\rm{~distinct}}} x^{\sigma(\lambda)}
  = \sum_{\sigma \in S_n; \sigma(\lambda) {\rm{~distinct}}} s_{\sigma(\lambda)}
\end{equation}
\end{formula}
\noindent{\bf{Proof.}}\quad 
The formula holds if and only if 
\begin{equation}
\sum_{\sigma \in S_n} x^{\sigma(\lambda)}
  = \sum_{\sigma \in S_n} s_{\sigma(\lambda)} \, ,
\end{equation}
since summing over all elements of $S_n$ adds the same symmetry factor
on each side of the equation.  
If we insert
\begin{equation}
s_{\sigma(\lambda)} = \frac{\sum_{w \in S_n} \epsilon(w) \,
x^{w(\sigma(\lambda)+\delta)}}
{\sum_{w \in S_n} \epsilon(w) \, x^{w(\delta)}} \, ,
\end{equation}
where $\epsilon(w)$ is the sign of the permutation $w$, 
we get
\begin{equation}
\sum_{\sigma,w \in S_n} \epsilon(w) \, x^{\sigma(\lambda)+w(\delta)} = 
\sum_{\sigma,w \in S_n} \epsilon(w) \, x^{w(\sigma(\lambda)+\delta)}  \, .
\end{equation}
Now, letting $w\sigma=\sigma'$ in the RHS of this equation, we obtain
\begin{equation}
\sum_{\sigma,w \in S_n} \epsilon(w) x^{\sigma(\lambda)+w(\delta)} = 
\sum_{\sigma',w \in S_n} \epsilon(w) x^{\sigma'(\lambda)+w(\delta)}  \, ,
\end{equation}
which proves the formula.
\hfill $\square$


\begin{thebibliography}{33}   
\bibitem{[Ad]}
A. M. Garsia, \emph{Lecture notes: Rodriguez formulas and orthogonality
for Schur, Hall-Littlewood, and Macdonald polynomials}, (1997).
\bibitem{[Ga]}
A. M. Garsia and M. Haiman, \emph{A graded representation module
for Macdonald's polynomials},
Proc. Natl. Acad. Sci. USA {\bf 90} (1993) 3607-3610.
\bibitem{[Ha]}
M. Haiman, \emph{Hilbert schemes, polygraphs, and the Macdonald
positivity conjecture}, J. Amer. Math Soc. {\bf{14}} (2001), 941--1006.
\bibitem{[Ji]} N. Jing, \emph{Vertex operators and Hall-Littlewood symmetric
functions}, Adv. Math. {\bf 87} (1991), 226--248.
\bibitem{[LLM]} L. Lapointe, A. Lascoux and J. Morse, 
\emph{Tableau atoms and a new Macdonald positivity conjecture}, to be published
in Duke Math. J., math.QA/0008073.
\bibitem{[LM2]} L. Lapointe and J. Morse, \emph{Schur function analogs for 
a filtration of the symmetric function space}, math.CO/0111192.
\bibitem{[Lint]} A. Lascoux, \emph{Interpolation}, 
http://www-igm.univ-mlv.fr/{\~{}}lascoux.
\bibitem{[Ma]} I.~G. Macdonald, {Symmetric Functions and Hall
    Polynomials}, 2nd edition, Clarendon Press, Oxford, 1995.
\bibitem{[Shi]}
A. Schilling and S. Warnaar,
\emph{Inhomogeneous lattice paths, generalized
Kostka-Foulkes polynomials, and
$A_{n-1}$-supernomials}, Comm. Math. Phys. {\bf 202} (1999), 359--401.
\bibitem{[S2]} M. Shimozono, \emph{A cyclage poset structure for
Littlewood-Richardson tableaux}, European J.
Combin. {\bf 22} (2001), 365--393.
\bibitem{[S1]} M. Shimozono and J. Weyman, \emph{Graded characters of modules supported in the closure of a nilpotent conjugacy class}, European J.
Combin. {\bf 21} (2000), 257--288.
\bibitem{[SZ]} M. Shimozono and M. Zabrocki, \emph{Hall-Littlewood vertex operators and generalized Kostka polynomials}, Adv. Math. {\bf 158} (2001), 66--85.
\bibitem{[V]} {S. Veigneau}, {\it ACE, an Algebraic
                      Combinatorics Environment for the computer
                      algebra system MAPLE\/},
{\it Version 3.0\/},1998, http://phalanstere.univ-mlv.fr/{\~{}}ace/.
 \end{thebibliography}
\end{document}